\documentclass[12pt]{amsart}

\usepackage{begnac}
\usepackage{stmaryrd}

\DeclareMathOperator{\BiPert}{BiPert}

\title{Continuous first order logic for unbounded metric structures}

\author{Itaï \textsc{Ben Yaacov}}

\address{Itaï \textsc{Ben Yaacov} \\
  Université Claude Bernard -- Lyon 1 \\
  Institut Camille Jordan \\
  43 boulevard du 11 novembre 1918 \\
  69622 Villeurbanne Cedex \\
  France}

\urladdr{\url{http://math.univ-lyon1.fr/~begnac/}}

\thanks{Research partially supported by NSF grant DMS-0500172,
  ANR chaire d'excellence junior THEMODMET (ANR-06-CEXC-007) and
  by Marie Curie research network ModNet.}
\thanks{The author would like to thank
  C.\ Ward Henson for many helpful discussions and comments.}

\svnInfo $Id: Unbdd.tex 866 2009-03-27 13:05:28Z begnac $
\thanks{\textit{Revision} {\svnInfoRevision} \textit{of} \today}


\begin{document}

\begin{abstract}
  We present an adaptation of continuous first order logic to
  unbounded metric structures.
  This has the advantage of being closer in spirit to
  C.\ Ward Henson's logic for Banach
  space structures than the unit ball approach
  (which has been the common approach so far to Banach space
  structures in continuous logic),
  as well as of applying in situations where the unit ball approach
  does not apply (i.e., when the unit ball is not a definable set).

  We also introduce the process of single point \emph{emboundment}
  (closely related to the topological single point compactification),
  allowing to bring unbounded structures back into the setting of
  bounded continuous first order logic.

  Together with results from \cite{BenYaacov:Perturbations} regarding
  perturbations of bounded metric structures,
  we prove a Ryll-Nardzewski style characterisation of theories of
  Banach spaces which are separably categorical up to small
  perturbation of the norm.
  This last result is motivated by an unpublished result of Henson.
\end{abstract}

\maketitle

\section*{Introduction}

Continuous first order logic is an extension of classical first
order logic, introduced in \cite{BenYaacov-Usvyatsov:CFO}
as a model theoretic formalism for metric structures.
It is convenient to consider that continuous logic also extends
C.\ Ward Henson's logic for Banach space structures
(see for example \cite{Henson-Iovino:Ultraproducts}),
even though this statement is obviously false:
continuous first order logic deals exclusively with
\emph{bounded} metric structures,
immediately excluding Banach spaces from the picture.
This is a technical hurdle which is relatively easy to overcome.
What one usually does
(e.g., in \cite[Example~4.5]{BenYaacov-Usvyatsov:CFO}
and the discussion that follows it) is decompose a
Banach space into a multi-sorted structure, with one sort
for, say, each closed ball of radius $n \in \bN$.
One may further rescale all such sorts into the sort of
the unit ball, which therefore suffices as a single sorted structure.
The passage between Banach space structures in Henson's logic
and unit ball structures in continuous logic preserves such
notions as elementary classes, elementary extensions,
type-definability of subsets of the unit ball, etc.
This approach has allowed so far to translate almost every
model theoretic question regarding Banach space structures
to continuous logic.

The unit ball approach suffers nonetheless from several drawbacks.
One drawback, which served as our original motivation,
comes to light in the context of perturbations of metric
structures introduced in \cite{BenYaacov:Perturbations}.
Specifically, we wish to consider the notion of perturbation of
the norm of a Banach space arising from the Banach-Mazur distance.
However, any linear isomorphism of Banach spaces which respects the
unit ball is necessarily isometric, precluding any possibility of a
non trivial Banach-Mazur perturbation.
Another drawback of the unit ball approach, also remedied by the
tools introduced in the present paper,
is that in some unbounded metric structures the unit ball is not
a definable set (even though it is always type-definable),
so naming it as a sort (and quantifying over it)
adds undesired structure.
For example, this is the case with complete normed fields
(i.e., of fields equipped with a complete non trivial multiplicative
valuation in $\bR$),
considered in detail in \cite{BenYaacov:NormedFields}.

In the present paper we replace the unit ball approach with the
formalism of \emph{unbounded continuous first order logic},
directly applicable to unbounded metric structures
and in particular to Banach space structures.
Using some technical definitions introduced in \fref{sec:Gauged},
the syntax and semantics of unbounded logic are defined in
\fref{sec:UnboundedLogic}.
In \fref{sec:Compactness} we prove Łoś's Theorem for unbounded logic,
and deduce from it a Compactness Theorem inside bounded sets.
It follows that the type space of an unbounded theory is locally
compact.
In \fref{sec:Henson} we show that unbounded continuous first order
logic has the same expressive power as Henson's logic of
positive bounded formulae.

In order to be able to apply to unbounded structures
tools which are already developed in the context of
standard (i.e., bounded) continuous logic, we introduce in
\fref{sec:Emboundment} the process of \emph{emboundment}.
Trough the addition of a single point at infinity, to each unbounded
metric structure we associate a bounded one, to which established
tools apply.
This method is used in \fref{sec:Perturbation} to adapt
the framework of perturbations, developed in
\cite{BenYaacov:Perturbations} for bounded structures,
to unbounded ones.
In particular, \fref{thm:UnbddPertRN} asserts that
the Ryll-Nardzewski style characterisation of
$\aleph_0$-categoricity up to perturbation
\cite[Theorem~3.5]{BenYaacov:Perturbations}
holds for unbounded metric structures as well.

As an application, we prove in \fref{sec:HensonCategoricity} a
Ryll-Nardzewski style characterisation of theories of Banach spaces
which are $\aleph_0$-categorical up to arbitrarily small perturbation
of the norm.
This result is motivated by an unpublished result of Henson,
whom we thank for the permission to include it in the
present paper.

\medskip

Notation is mostly standard.
We use $a$, $b$, $c$, \ldots to denote members of structures,
and use $x$, $y$, $z$, \ldots to denote variables.
Bar notation is used for (usually finite) tuples, and uppercase
letters are used for sets.
We also write $\bar a \in A$ to say that $\bar a$ is a tuple
consisting of members of $A$, i.e., $\bar a \in A^n$ where
$n = |\bar a|$.
When $T$ is an $\cL$-theory (whether bounded or unbounded)
we always assume that $T$ is closed under logical consequences.
In particular, $|T| = |\cL| + \aleph_0$ and
$T$ is countable if and only if $\cL$ is.
We shall assume familiarity with (bounded) continuous first order
logic, as developed in \cite{BenYaacov-Usvyatsov:CFO}.
For the parts dealing with perturbations, familiarity with
\cite{BenYaacov:Perturbations} is assumed as well.
For a general survey of the model theory of metric structures
we refer the reader to
\cite{BenYaacov-Berenstein-Henson-Usvyatsov:NewtonMS}.

\section{Gauged spaces}
\label{sec:Gauged}

We would like to allow unbounded structures, while at the same time
keeping some control over the behaviour of bounded parts thereof.
The ``bounded parts'' of a structure are given by means of
a \emph{gauge}.

\begin{dfn}
  \label{dfn:NuSpace}
  Let $(X,d)$ be a metric space,
  $\nu\colon X \to \bR$ any function.
  We define $X^{\nu\leq r} = \{x \in X\colon \nu(x) \leq r\}$ and similarly 
  $X^{\nu\geq r}$, $X^{\nu<r}$, etc.
  \begin{enumerate}
  \item We call
    $X^{\nu\leq r}$ and $X^{\nu<r}$ the \emph{closed} and \emph{open $\nu$-balls}
    of radius $r$ in $X$, respectively.
  \item We say that $\nu$ is a \emph{gauge} on $(X,d)$,
    and call the triplet
    $(X,d,\nu)$ a \emph{($\nu$-)gauged space} if
    $\nu$ is $1$-Lipschitz in $d$ and
    every $\nu$-ball (of finite radius) is bounded in $d$.
  \end{enumerate}
  Note that this implies
  that the bounded subsets of $(X,d)$ are precisely those contained
  in some $\nu$-ball.
\end{dfn}

\begin{rmk}
  We could have given a somewhat more general definition, replacing
  the $1$-Lipschitz condition with the weaker condition that
  the gauge
  $\nu$ should be bounded and uniformly continuous on every bounded set.
  This does not cause any real loss of generality, since in that case
  we could define
  \begin{gather*}
    d'(x,y) = d(x,y) + |\nu(x) - \nu(y)|.
  \end{gather*}
  Then $\nu$ is $1$-Lipschitz with respect to $d'$,
  and the two metrics $d$ and $d'$ are uniformly equivalent and
  induce the same notion of a bounded set.
\end{rmk}

\begin{dfn}
  Recall that a \emph{(uniform) continuity modulus} is a
  left-continuous increasing function
  $\delta\colon (0,\infty) \to (0,\infty)$
  (i.e.,
  $\delta(\varepsilon)
  = \sup_{\varepsilon<\varepsilon'} \delta(\varepsilon')$).

  We say that a mapping
  $f\colon (X,d_X,\nu_X) \to (Y,d_Y,\nu_Y)$ between two gauged spaces
  \emph{respects $\delta$ under $\nu$} if for all
  $\varepsilon > 0$:
  \begin{gather*}
    \tag{UC$_\nu$}\label{eq:UCnu}
    \begin{aligned}
      & \nu_X(x),\nu_X(y) < \inv{\varepsilon}, \\
      & d_X(x,y) < \delta(\varepsilon)
    \end{aligned}
    \quad \Longrightarrow \quad
    \begin{aligned}
      & d_Y(f(x),f(y)) \leq \varepsilon, \\
      & \nu_Y(x) \leq \inv{\delta(\varepsilon)}.
    \end{aligned}
  \end{gather*}
  We say that $f$ is \emph{uniformly continuous under $\nu$}
  if it respects some $\delta$ under $\nu$.
\end{dfn}

While respecting a given $\delta$ under $\nu$ depends on the choice
of $\nu$, the fact that some $\delta$ is respected under $\nu$ does
not.
\begin{lem}
  Let $X$ and $Y$ be gauged spaces,
  $f\colon X \to Y$ a mapping.
  Then
  \begin{enumerate}
  \item
    Let $\delta\colon (0,\infty) \to (0,\infty)$ be any mapping, and
    assume that $f$ respects $\delta$ under $\nu$ in the sense of
    \fref{eq:UCnu}.
    Define
    $\delta'(\varepsilon)
    = \varepsilon
    \wedge \sup_{0<\varepsilon'<\varepsilon} \delta(\varepsilon')$.
    Then $\delta' \leq \id$ is a continuity modulus and
    $f$ respects $\delta'$ under $\nu$ as well.
    (If we used $\sup$ alone we could obtain infinite values, whence
    the need for truncation at $\varepsilon$.)
  \item
    A mapping $f\colon X \to Y$ between gauged spaces
    is uniformly continuous under $\nu$ if and only if it restriction
    to every bounded set is
    uniformly continuous and bounded.
  \end{enumerate}
\end{lem}
\begin{proof}
  Easy.
\end{proof}

\begin{dfn}
  \label{dfn:GaugeCartesianProduct}
  A Cartesian product of gauged metric spaces
  $X = \prod_{i<n} X_i$ is equipped with a gauged metric structure as
  follows:
  \begin{gather}
    d(\bar x,\bar y) = \bigvee_{i<n} d(x_i,y_i),
    \qquad
    \nu(\bar x) = \bigvee_{i<n} \nu(x_i).
  \end{gather}
  In particular, if $n = 0$ then $X = \{*\}$
  and $d(*,*) = \nu(*) = 0$.

  We also identify $\bR^+$ with the gauged space
  $(\bR^+,|x-y|,|x|)$.
\end{dfn}

\begin{lem}
  \label{lem:GaugedUniformContinuity}
  Let $X$, $Y$, and so on, denote gauged spaces.
  \begin{enumerate}
  \item
    The projection mapping $X \times Y \to X$ respects the identity
    uniformly under $\nu$.
  \item
    Let $f_i\colon X \to Y_i$, $i < n$,
    be mappings between gauged spaces,
    each respecting $\delta_{f_i}$ under $\nu$.
    Then
    $\bar f\colon X \to \prod Y_i$
    respects the continuity modulus
    $\delta_{\bar f} = \bigwedge_{i<n} \delta_{f_i}$ under $\nu$.
    In addition, if $\delta_{f_i} \leq \id$ for all (indeed, for some)
    $i < n$ then $\delta_{\bar f} \leq \id$ as well.
  \item Let $X$, $Y$ and $Z$ be gauged spaces.
    Assume that $f\colon X \to Y$ and $g\colon Y \to Z$ respect
    continuity moduli $\delta_f$ and $\delta_g$, respectively,
    under $\nu$.
    Assume moreover that $\delta_f,\delta_g \leq \id$.
    Then $h = g \circ f\colon X \to Z$ respects the
    continuity modulus
    $\delta_h = \delta_f \circ \delta_g \circ \delta_f$
    under $\nu$.
    In particular, $\delta_h \leq \id$ is a continuity modulus.
  \item Let $X$ and $Y$ be gauged spaces,
    and let $f\colon X \times Y \to \bR^+$ and
    $g\colon Y \to \bR^+$ mappings
    which respect $\delta_f$ and $\delta_g$ under $\nu$, respectively.
    Assume also that $f$ is \emph{eventually equal} to $g$,
    namely that there exists a constant $C$ such that
    $f(x,y) = g(y)$ whenever $\nu(x) \geq C$.
    Define
    \begin{gather*}
      \begin{aligned}
        h_1(y) & = \sup_{x\in X} f(x,y), \qquad &
        h_1'(y) & = g(y) \vee \sup_{x\in X} f(x,y),
        \\
        h_2(y) & = \inf_{x\in X} f(x,y), &
        h_2'(y) & = g(y) \wedge \inf_{x\in X} f(x,y),
      \end{aligned}
      \\
      \delta_h(\varepsilon)
      = \delta_g(\varepsilon)
      \wedge \delta_f(\varepsilon \wedge \inv C).
    \end{gather*}
    Then $h_i, h_i'\colon Y \to \bR^+$ are well defined
    (i.e., the supremum is always finite)
    and respect $\delta_h$ under $\nu$.
    Moreover, if either $\delta_f \leq \id$ or $\delta_g \leq \id$
    then $\delta_h \leq \id$.
  \end{enumerate}
\end{lem}
\begin{proof}
  The first two items are easy.

  For the third item we only prove that $\delta_h$ is respected under
  $\nu$.
  Indeed, let $\varepsilon > 0$, $x,y \in X$,
  and assume that
  $\nu(x),\nu(y) < \inv{\varepsilon}$
  and $d(x,y) < \delta_h(\varepsilon)$.
  By the left continuity assumption there are
  $s,t$ such that:
  $d(x,y)
  < \delta_f(s)$,
  $s < \delta_g(t)$,
  $t < \delta_f(\varepsilon)$.
  In particular $s < t < \varepsilon$.
  Using our hypotheses we obtain from top to bottom:
  \begin{align*}
    & d(x,y) < \delta_f(s),
    &&
    \nu(x),\nu(y) < \inv{\varepsilon} < \inv s,
    \\
    & d\bigl( f(x), f(y) \bigr) \leq s < \delta_g(t),
    &&
    \nu\circ f(x), \nu \circ f(y)
    \leq \inv{\delta_f(\varepsilon)} < \inv t,
    \\
    & d\bigl( h(x), h(y) \bigr) \leq t < \varepsilon,
    &&
    \nu\circ h(x), \nu \circ h(y) \leq \inv{\delta_g(t)} < \inv s.
  \end{align*}
  In addition, $s$ could have been chosen arbitrarily close to
  $\delta_g \circ \delta_f(\varepsilon)$
  whereby
  $\nu\circ h(x)
  \leq \inv{\delta_g \circ \delta_f(\varepsilon)}
  \leq \inv{\delta_f \circ \delta_g \circ \delta_f(\varepsilon)}
  = \inv{\delta_h(\varepsilon)}$,
  as desired.

  For the fourth item, the existence of $h_1$ and $h_1'$
  follows from the fact
  that for a fixed $y$, the function $x \mapsto f(x,y)$
  is bounded on bounded sets and eventually constant.
  We show that $h_1$ respects $\delta_h$ under $\nu$, a similar
  argument applies to the other functions.
  Let $y_1,y_2 \in Y$,
  and assume that
  $\nu(y_i) < \inv{\varepsilon}$,
  $d(y_1,y_2) < \delta_h(\varepsilon)$.
  Let $r = \varepsilon \wedge \inv C$,
  so $\nu(y_i) < \inv r$ and
  $d(y_1,y_2) < \delta_g(\varepsilon) \wedge \delta_f(r)$.
  We may choose a point $x \in X$ such that
  $f(x,y_1)$ is arbitrarily close to $h_1(y_1)$.
  There are two cases to consider:
  \begin{align*}
    & \textbf{I. } \nu(x) \geq C
    &&
    \begin{aligned}
      & f(x,y_1) = g(y_1)
      \leq \inv{\delta_g(\varepsilon)}
      \leq \inv{\delta_h(\varepsilon)}, \\
      &
      |f(x,y_1)-f(x,y_2)| = |g(y_1)-g(y_2)| \leq \varepsilon,
    \end{aligned}
    \\
    &
    \\
    & \textbf{II. } \nu(x) < C \leq \inv r
    &&
    \begin{aligned}
      & f(x,y_1)
      \leq \inv{\delta_f(r)}
      \leq \inv{\delta_h(\varepsilon)},
      \\
      & |f(x,y_1)-f(x,y_2)| \leq r \leq \varepsilon.
    \end{aligned}
  \end{align*}
  Either way we obtain that
  $h_1(y_1) \leq \inv{\delta_h(\varepsilon)}$
  and that
  $h_1(y_1) \leq h_1(y_2) + \varepsilon$,
  which is enough.

  If there exists $x \in X$ such that $\nu(x) \geq C$ then
  $h_1 = h_1'$.
  If not then when dealing with $h_1'$ we need to consider the
  possibility that $h_1'(y_1) = g(y_1)$,
  which is treated identically to case \textbf{I}.
  The functions $h_2$ and $h_2'$ are treated analogously.
\end{proof}

\section{Unbounded continuous logic}
\label{sec:UnboundedLogic}

We turn to define a $\bR^+$-valued variant of continuous logic which
can accommodate unbounded metric structures.
We shall refer to this logic as \emph{unbounded continuous logic}.
The $[0,1]$-valued (or, more generally, bounded) continuous logic
defined in \cite{BenYaacov-Usvyatsov:CFO} will be referred to here as
\emph{standard} or \emph{bounded}.

\begin{dfn}
  An \emph{unbounded continuous signature} $\cL$
  consists of the following data:
  \begin{enumerate}
  \item A set of \emph{relation (or predicate) symbols} and of
    \emph{function symbols}, each equipped with its \emph{arity}
    (zero-ary function symbols are also called
    \emph{constant symbols}).
  \item
    For each $n$-ary symbol $s$,
    a continuity modulus
    $\delta_s\colon (0,\infty) \to (0,\infty)$.
  \item 
    For each sort $S$, a distinguished binary predicate symbol
    $d_S$ called the \emph{distance symbol},
    as well as a distinguished
    unary predicate symbol $\nu_S$ called the \emph{gauge symbol}.
    The subscript $S$ is usually omitted.
  \end{enumerate}
  We usually write down a signature merely by listing its
  non distinguished symbols.
\end{dfn}

\begin{dfn}
  Let $\cL$ be an unbounded signature, and for the sake of simplicity
  let us assume it is single-sorted.
  An \emph{(unbounded) $\cL$-structure} is a complete metric gauged
  space $(M,d,\nu) = (M,d^M,\nu^M)$, possibly empty,
  equipped with interpretation of the symbols:
  \begin{enumerate}
  \item
    The interpretation of an $n$-ary function symbol $f$ is a
    mapping $f^M \colon M^n \to M$
    which respects $\delta_f$ under $\nu$.
  \item
    The interpretation of an $n$-ary predicate symbol $P$ is a
    mapping $P^M \colon M^n \to \bR^+$
    which respects $\delta_P$ under $\nu$.
  \end{enumerate}
  For this purpose we view $M^n$ with a gauged space
  $(M^n,d,\nu)$ as per \fref{dfn:GaugeCartesianProduct}.
  Similarly, $\bR^+$ admits a standard gauge structure
  $(\bR^+,d,\id)$.
\end{dfn}

Thus, restricted to a $\nu$-ball, everything is bounded and uniformly
continuous as in bounded continuous logic, and closed $\nu$-balls are
metrically closed and therefore complete.

\begin{rmk}
  If the language contains a constant symbol $0$
  then the formula $\nu'(x) = d(x,0)$ can act as an alternative gauge.
  Indeed, if $r \in \bR^+$ then
  $M^{\nu'\leq r} \subseteq M^{\nu\leq r+\nu(0)}$, since $\nu$ is
  $1$-Lipschitz, and conversely
  $M^{\nu\leq r} \subseteq M^{\nu'\leq\beta_d(r,\nu(0))}$ by
  definition of an unbounded structure.
  Thus we can pass between $\nu$-balls and $\nu'$-balls in a way which
  depends only on $\cL$.

  In most cases, $\nu$ will indeed be equal to $d(x,0)$.
\end{rmk}

A \emph{standard} continuity modulus for an $n$-ary symbol,
when $n > 0$, is the function $x \mapsto \frac{x}{n}$.
If a symbol $s$ is $1$-Lipschitz is each argument
and $\nu(x) = d(x,0)$ then $s$ indeed respects the standard
continuity modulus under $\nu$.
For a zero-ary symbol the standard continuity modulus is the
identity.

\begin{exm}
  \label{exm:BoundedStructreUnbounded}
  Let $\cL$ be a standard (i.e., $[0,1]$-valued) continuous signature
  as defined in \cite{BenYaacov-Usvyatsov:CFO}.
  In that case we chose to equip each $n$-ary symbol $s$ with
  individual continuity moduli $\delta_{s,i}$, $i < n$,
  one for each argument.
  Let $\cL'$ be the unbounded signature obtained from $\cL$ by adding
  a gauge symbol $\nu_S$ for each sort $S$, and by setting
  $\delta_s(x) = 1 \wedge \bigwedge_i \delta_{s,i}(\frac{x}{n})$
  (for zero-ary $s$ let $\delta_s = 1$).
  Then every $\cL$-structure $M$ can be naturally viewed as an
  unbounded $\cL'$-structure by interpreting all gauges as the
  constant $0$.
  If $\cL$ admits a constant symbol $0$
  then interpreting $\nu(x) = d(x,0)$ works as well.
\end{exm}

\begin{exm}[Banach spaces]
  \label{exm:BanachSig}
  We would like to view Banach spaces as unbounded structures.
  Let $\cL = \{0,+,m_r\colon r \in \bQ\}$, where $m_r$ is
  unary scalar multiplication by $r$.
  We view $\|x\|$ as shorthand for $d(x,0)$, and take it to be the
  gauge.
  Let $\delta_{m_r}(x) = |r|x$
  and let all other continuity moduli be standard.
  Then every real Banach space is naturally an (unbounded)
  $\cL$-structure.

  This can be extended to additional structure on the Banach space.
  For example a complex Banach space also has a function symbol for
  multiplication by $i$, while a Banach lattice is given by binary
  function symbols $\vee$, $\wedge$
  (again with standard continuity moduli).
\end{exm}

\begin{exm}[Naming constants]
  \label{exm:NamingConstants}
  Let $\cL$ be an unbounded signature, $M$ and $\cL$-structure.
  Let $A \subseteq M$.
  We define $\cL(A)$ as $\cL \cup A$, where each $a \in A$ is viewed
  as a new constant symbol.
  We equip each symbol $a$ with the uniform continuity modulus
  $\delta_a = \id \wedge \frac{1}{\nu(a)}$.
  Then for every $\varepsilon > 0$ we have
  $\nu(a) \leq \frac{1}{\delta_a(\varepsilon)}$, and we may render $M$
  an $\cL(A)$-structure by interpreting
  $a^M = a$.
\end{exm}

We now define the syntax of continuous logic.
A term is defined, as usual, as either being a variable or
a composition of a function symbol with simpler terms.
Similarly, an atomic formula is a composition of a predicate symbol
with terms.
Connectives are continuous functions from
$(\bR^+)^n$ to $\bR^+$, or any convenient family of such functions
which is dense in the compact-open topology, i.e., in the topology of
uniform convergence on every compact set.
We shall use the system
$\{1,x \dotminus y,x+y,x/2\}$ which  generates such a dense set
through composition.
While alternative systems may be legitimate,
we shall always require the presence of $1$ and $\dotminus$ in what
follows.
We point out that as functions from $(\bR^+)^n \to \bR^+$,
all the basic connectives we chose respect their respective
standard continuity moduli (according to their arity).
As one may expect, every combination of formulae by connectives is a
formula.

On the other hand, care is needed when defining quantified formulae.
First, $\sup_x \varphi$ could be infinite.
Second, even if $\varphi$ is bounded, we still need a uniform rate of
convergence for
$\sup_{\nu(x) < C} \varphi \to \sup_x \varphi$ as $C \to \infty$,
or else we may run into trouble
with compactness as well as with uniform continuity under $\nu$.
In Henson's logic of positive bounded formulae
\cite{Henson-Iovino:Ultraproducts}, where
the truth values are True/False, one gets around this by restricting
quantifiers to bounded balls
(and then again, one needs to play around with the radii of the balls
when considering approximations).
If we tried to do the same thing with continuous quantifiers we
could again run into trouble if, say,
$\sup_{\nu(x)<r} \varphi < \sup_{\nu(x)\leq r} \varphi$.
We shall follow a different path, looking
for the simplest \emph{syntactic} conditions on a
formula $\varphi$ that ensure that $\inf_x \varphi$ and $\sup_x \varphi$ are
\emph{semantically} legitimate.
This approach will allow us nonetheless to
recover approximate versions of bounded quantifiers later on.

\begin{dfn}
  \label{dfn:UnbddFormula}
  We define formulae by induction, and at the same time we
  define whether a formula is syntactically
  eventually constant in a variable $x$ and/or bounded.
  \begin{itemize}
  \item Atomic formulae are defined as above.
    \\
    -- If $\varphi$ is atomic and $x$ does not appear in $\varphi$,
    then $\varphi$ is eventually constant in $x$.
    \\
    -- No atomic formula is bounded.
  \item A combination of formulae by connectives is a formula.
    \\
    -- If all the components are bounded
    (respectively, eventually constant in $x$)
    then so is the compound formula.
    \\
    -- If $\varphi$ is bounded then $\varphi \dotminus \psi$ is
    bounded for any $\psi$ and $\varphi \dotminus \nu(x)$ is
    eventually constant in $x$.
  \item If $\varphi$ is eventually constant in $x$ then
    $\inf_x \varphi$ and $\sup_x \varphi$ are formulae
    (but not otherwise).
    \\
    -- If $\varphi$ is bounded
    (respectively, eventually constant in a variable $y$)
    then so are $\sup_x \varphi$ and $\inf_x \varphi$.
    In particular, $\sup_x \varphi$ and $\inf_x \varphi$
    are eventually constant in $x$.
  \end{itemize}
\end{dfn}

Notice that the formula $1$, being a combination of no formulae,
is bounded and eventually constant in every variable.
Similarly, every dyadic number
$r = \frac{k}{2^m} = (\half)^m(1 + \cdots + 1)$
can be viewed as a formula, and is syntactically bounded and constant
as such.
It follows that for every formula $\varphi$,
the formula
$\varphi\wedge r = r\dotminus (r \dotminus \varphi)$
is syntactically bounded.

The qualitative syntactic properties of boundedness and eventual
constancy can be translated to quantitative information.
\begin{dfn}
  \label{dfn:QuantitativeBounds}
  For every syntactically bounded formula $\varphi$ we extract a
  syntactic bound $B_\varphi$ as follows:
  \begin{align*}
    & \varphi = \theta(\bar \psi), \text{ $\psi_i$ bounded:}
    && B_\varphi = \sup_{\bar x \in \prod [0,B_{\psi_i}]} \theta(\bar x),
    \\
    & \varphi = \psi \dotminus \chi,
    \text{ or } \sup_x \psi,
    \text{ or } \inf_x \psi, \text{ $\psi$ bounded:}
    && B_\varphi = B_\psi.
  \end{align*}
  Notice that no ambiguity arises for $\psi \dotminus \chi$ when both
  $\psi$ and $\chi$ are syntactically bounded.

  Similarly, for a formula $\varphi(x,\bar y)$ which is syntactically
  eventually constant in $x$ we extract a syntactic constancy
  threshold $C_{\varphi,x} \in \bR^+$ and a formula
  $\varphi(\infty,\bar y)$, whose free variables lie among $\bar y$,
  and which is intended to agree with $\varphi(x,\bar y)$
  once $\nu(x) \geq C_{\varphi,x}$.
  \begin{align*}
    & \text{$x$ not free in $\varphi$:}
    && \varphi(\infty,\bar y) = \varphi,
    && C_{\varphi,x} = 0,
    \\
    & \varphi = \theta(\bar \psi), \text{ $\psi_i$ e.c.\ in $x$:}
    && \varphi(\infty,\bar y) = \theta\bigl( \psi(\infty,\bar y) \bigr),
    && C_{\varphi,x} = \bigvee C_{\psi_i,x},
    \\
    & \varphi = \psi \dotminus \nu(x), \text{ $\psi$ bounded:}
    && \varphi(\infty,\bar y) = 0,
    && C_{\varphi,x} = B_\psi,
    \\
    & \varphi = \sup_z \psi, \text{ $\psi$ e.c.\ in $x \neq z$:}
    && \varphi(\infty,\bar y) = \sup_z \psi(\infty,z,\bar y),
    && C_{\varphi,x} = C_{\psi,x}.
  \end{align*}
  Again, when cases overlap the definitions agree.
\end{dfn}

The definition of the semantics can be somewhat delicate.
The model for the definition is an unbounded structure $M$ in which
elements of arbitrarily high gauge exist (e.g., a non trivial Banach
space).
In this case the definition is entirely straightforward,
namely
\begin{align*}
  & (f\bar \tau)^M(\bar a) = f^M \circ \bar \tau^M(\bar a) \in M,
  &&
  (P\bar \tau)^M(\bar a) = P^M \circ \bar \tau^M(\bar a) \in \bR^+,
  \\
  & \theta(\bar \varphi)^M(\bar a)
  = \theta\bigl( \bar \varphi^M(\bar a) \bigr)
  && \qquad\qquad \text{(where $\theta$ is a connective)},
  \\
  \tag{Q}\label{eq:QuantifierWithoutInfinity}
  & \bigl( \inf_x \varphi(x,\bar a) \bigr)^M
  = \inf_{b\in M} \varphi^M(b,\bar a),
  &&
  \text{idem for } \sup.
\end{align*}
Let us state some properties of this model situation, for the time
being without proof.
First, the interpretation of every term and formula is uniformly
continuous under $\nu$
(essentially by \fref{lem:GaugedUniformContinuity}).
Second, if $\varphi$ is syntactically bounded then it is bounded by
$B_\varphi$.
Third, if $\varphi(x,\bar y)$ is syntactically eventually constant
in $x$ then $\varphi(x,\bar y) = \varphi(\infty,\bar y)$
whenever $\nu(x) \geq C_{\varphi,x}$.
In this case, $\varphi(x,\bar a)$ is bounded for every $\bar a$, so
the interpretation of the quantifiers makes sense
and the following holds:
\begin{align*}
  \tag{Q$_\infty$}\label{eq:QuantifierWithInfinity}
  & \bigl( \inf_x \varphi(x,\bar a) \bigr)^M
  = \inf_{b\in M\cup \{\infty\}} \varphi^M(b,\bar a),
  &&
  \text{idem for } \sup.
\end{align*}
However, we must also take into account structures
in which elements of arbitrarily high gauge need not
exist.
In order for ultra-products to behave reasonably,
i.e., in order for Łoś's Theorem to hold, the definition of quantifier
semantics in the general case must follow
\fref{eq:QuantifierWithInfinity} and \emph{not}
\fref{eq:QuantifierWithoutInfinity}.
This is illustrated in \fref{rmk:UnboundedPathologies} below.

\begin{dfn}
  Let $M$ be an $\cL$ structure.
  Then terms, atomic formulae and connectives are interpreted
  naturally, by composition.
  Quantifiers are interpreted according to
  \fref{eq:QuantifierWithInfinity},
  where $\varphi(\infty,\bar a)$ is understood as per
  \fref{dfn:QuantitativeBounds}.
\end{dfn}

\begin{thm}
  Let $M$ be an $\cL$-structure.
  Then:
  \begin{enumerate}
  \item All formulae are interpreted as $\bR^+$-valued functions on
    Cartesian powers of $M$.
    In particular, in the interpretation of quantified formulae in $M$
    all the suprema are finite.
  \item Every term $\tau$ and every formula $\varphi$ are uniformly
    continuous under $\nu$.
  \item If a formula $\varphi$ is syntactically constant then
    $\varphi^M(\bar a) \leq B_\varphi$ for all $\bar a \in M$.
  \item If a formula $\varphi(x,\bar y)$ is
    syntactically eventually constant in $x$ then
    $\varphi(b,\bar a) = \varphi(\infty,\bar a)$
    whenever $\nu(b) \geq C_{\varphi,x}$.
  \end{enumerate}
\end{thm}
\begin{proof}
  We prove this by induction on the complexity of terms and formulae.
  We observe that if $\varphi(x,\bar y)$ is syntactically eventually
  constant in $x$ then $\varphi(\infty,\bar y)$ is of lesser or equal
  complexity.
  Thus, when treating $\sup_x \varphi$ and $\inf_x \varphi$,
  we may use the induction hypotheses both for $\varphi(x,\bar y)$
  and for $\varphi(\infty,\bar y)$.
  We may assume that all the continuity moduli of symbols lie below
  the identity, and construct as we go continuity moduli
  below the identity for each term and formula.

  The induction step itself now follows immediately from the
  definitions, the induction hypotheses and
  \fref{lem:GaugedUniformContinuity}.
\end{proof}

We leave it as an exercise to the reader to check that
with our choice of connectives, every formula is equivalent
to one in prenex normal form
(one needs to make sure in particular that the natural transformations
towards a prenex form do not violate the restrictions on
quantification imposed in \fref{dfn:UnbddFormula}).

It will be convenient later on to have some
analogue of the restricted quantifier
$\sup_{\nu(x)\leq r} \varphi$ (which is not part of our language).
Let us assume that $\varphi$ is syntactically bounded
and let $k = \lceil B_\varphi \rceil$, namely the least integer
syntactic bound for $\varphi$.
We observe that for a dyadic $r$, the formula
$\varphi \dotminus (\nu(x)\dotminus r)$ is equivalent to
$(\varphi+(\nu(x)\wedge r)) \dotminus \nu(x)$ which is
syntactically bounded and eventually constant in $x$.
It follows that for every natural $m > 0$ the formula
$\varphi \dotminus m(\nu(x)\dotminus r)$ is equivalent to one which is
syntactically bounded and eventually constant in $x$.
Let $0 < r < r'$, and
find the least $m$ such that we can write
$r \leq s = \ell2^{-m} < (\ell+1)2^{-m} \leq r'$, and choose the least
possible $s$.
Define:
\begin{align*}
  & \varphi{\downarrow}^{x\leq r,r'}
  = \varphi \dotminus k2^m(\nu(x)\dotminus s),
  &&
  {\sup_x}^{r,r'} \varphi
  = \sup_x \varphi{\downarrow}^{x\leq r,r'},
  \\
  & \varphi{\uparrow}^{x\leq r,r'}
  = k \dotminus (k \dotminus \varphi){\downarrow}^{x\leq r,r'},
  &&
  {\inf_x}^{r,r'} \varphi
  = \inf_x \varphi{\uparrow}^{x\leq r,r'}.
\end{align*}
Both formulae on the left are
syntactically bounded and eventually constant in $x$,
so the expressions on the right are indeed formulae.
By construction we always have
$\varphi{\downarrow}^{x\leq r,r'} \leq \varphi$,
and in addition
$\varphi{\downarrow}^{x\leq r,r'} = \varphi$
when $\nu(x) \leq r$ and
$\varphi{\downarrow}^{x\leq r,r'} = 0$
when $\nu(x) \geq r'$.
Thus
$\sup_{\nu(x)\leq r} \varphi
\leq sup_x^{r,r'} \varphi
\leq \sup_{\nu(x) < r'} \varphi$.
Similarly,
$\inf_{\nu(x)\leq r} \varphi
\geq \inf_x^{r,r'} \varphi
\geq \inf_{\nu(x) < r'} \varphi$.

We may further extend these abbreviations to the case where $\varphi$ is not
syntactically bounded by truncating it at $1$, defining
$\varphi{\downarrow}^{x\leq r,r'}
= (\varphi\wedge1){\downarrow}^{x\leq r,r'}$
(and proceeding as above).
This will only be used in conditions of the form
$\sup_x^{r,r'} \inf_y^{s,s'} \ldots \varphi = 0$, whose satisfaction does not
depend on our particular choice of constant at which we truncate.

\section{Łoś's Theorem, compactness and theories}
\label{sec:Compactness}

Let $\cL$ be an unbounded signature,
$\{M_i\colon i \in I\}$ a family of $\cL$-structures and $\sU$
an ultra-filter on $I$.
Let $I_0 = \{i\colon M_i \neq \emptyset\}$.
If $I_0 \in \sU$ define
\begin{gather*}
  N_0
  = \left\{
    (a_i) \in \prod_{i\in I_0} M_i
    \colon \lim_{\sU} \nu^{M_i}(a_i) < \infty
  \right\},
\end{gather*}
otherwise $N_0 = \emptyset$.
Alternatively, one may introduce a new formal element
$\infty$ with $\nu(\infty) = +\infty$,
and define
\begin{gather*}
  N_0
  = \left\{
    (a_i) \in \prod_{i\in I} \bigl( M_i \cup \{\infty\} \bigr)
    \colon \lim_{\sU} \nu^{M_i}(a_i) < \infty
  \right\}.
\end{gather*}
Under this definition a member $(a_i) \in N_0$ can have few
(according to $\sU$) coordinates which are equal to $\infty$ and which
may be ignored in the definitions that follow.
Either approach leads to the same construction.

For a function symbol $f$ or predicate symbol $P$, and arguments
$(a_i),(b_i),\ldots \in N_0$, define:
\begin{gather*}
  f^{N_0}\bigl( (a_i),(b_i),\ldots \bigr)
  = (f^{M_i}(a_i,b_i,\ldots)),
  \\
  P^{N_0}\bigl( (a_i),(b_i),\ldots \bigr)
  = \lim_{\sU} P^{M_i}(a_i,b_i,\ldots).
\end{gather*}
Note that by definition of $N_0$, the values of $P^{M_i}(a_i,b_i,\ldots)$
are bounded on a large set of indexes, so
$\lim_{\sU} P^{M_i}(a_i,b_i,\ldots) \in \bR^+$.
It is now straightforward verification that $N_0$ is an
$\cL$-pre-structure, i.e., that it verifies all the properties of a
structure with the exception
that $d^{N_0}$ might be a pseudo-metric and needs not be complete.
Let $N = \widehat N_0$ be the associated $\cL$-structure,
obtained by dividing by the zero distance equivalence relation and
passing to the metric completion.
We call $N$ the
\emph{ultra-product of $\{M_i\colon i \in I\}$ modulo $\sU$},
denoted $\prod M_i/\sU$.
The image in $N$ of $(a_i) \in N_0$ will be denoted $[a_i]$.
(Compare with the construction of ultra-products of Banach spaces in
\cite{Henson-Iovino:Ultraproducts} and of bounded continuous
structures in \cite{BenYaacov-Usvyatsov:CFO}.)

\begin{thm}[Łoś's Theorem]
  For every formula $\varphi(\bar x)$ and $[a_i],[b_i],\ldots \in
  \prod M_i/\sU$:
  \begin{gather*}
    \varphi([a_i],[b_i],\ldots)^{\prod M_i/\sU}
    = \lim_\sU \varphi(a_i,b_i,\ldots)^{M_i}.
  \end{gather*}
\end{thm}
\begin{proof}
  Mostly as for bounded logic.
  The only significant difference is in the treatment of quantifiers,
  which we sketch below.

  If $\lim_\sU \inf_x \varphi(x,a_i,\ldots)^{M_i} < r$
  then there is a large set on which
  $\inf_x \varphi(x,a_i,\ldots)^{M_i} < r$
  and we can find witnesses $b_i$ there (possibly the formal infinity)
  such that
  $\varphi(b_i,a_i,\ldots)^{M_i} < r$.
  If $\nu(b_i) \leq C_{\varphi,x}$ on a large set
  then $\nu([b_i]) \leq C_{\varphi,x}$,
  so in particular $[b_i]$ belongs to the ultra-product and
  \begin{gather*}
    \inf_x \varphi(x,[a_i],\ldots)
    \leq \varphi([b_i],[a_i],\ldots)
    = \lim_{\sU} \varphi(b_i,a_i,\ldots) \leq r.
  \end{gather*}
  If, on the other hand,
  $b_i = \infty$ or $\nu(b_i) \geq C_{\varphi,x}$ on a large set
  then
  \begin{gather*}
    \inf_x \varphi(x,[a_i],\ldots)
    \leq \varphi(\infty,[a_i],\ldots)
    = \lim_{\sU} \varphi(\infty,a_i,\ldots)
    = \lim_{\sU} \varphi(b_i,a_i,\ldots) \leq r.
  \end{gather*}
  Conversely, assume that
  $\inf_x \varphi(x,[a_i],\ldots) < r$.
  Then again, either there is $[b_i]$ such that
  $\varphi([b_i],[a_i],\ldots) < r$
  or $\varphi(\infty,[a_i],\ldots) < r$,
  and in either case
  $\lim_{\sU} \inf_x \varphi(x,a_i,\ldots) \leq r$.
\end{proof}

\begin{rmk}
  \label{rmk:UnboundedPathologies}
  Łoś's Theorem might fail if our semantic interpretation did not take
  the value at infinity into account.
  For example, consider the sentence $\varphi = \inf_x (1 \dotminus \nu(x))$.
  Let $M_n$ be the structure consisting of two points,
  $\nu(a_{n,0}) = 0$, $\nu(a_{n,1}) = n$.
  Then the ultra-product contains a single point
  $a_0 = [a_n,0]$, $\nu(a_0) = 0$, and we would have
  $\varphi^{M_n} = 0$ for all $n \geq 1$
  and yet $\varphi^{\prod M_n/\sU} = 1$.

  Worse still, if $M_n$ consisted only of $a_{n,1}$ then
  $\prod M_i/\sU$ would be empty, making the naïve interpretation of
  quantifiers meaningless.
  An empty ultra-product 
  can also be obtained with unbounded structures, for example
  $M_n = E \setminus B(n)$ where $E$ is a Banach space and $B(n)$ is its open
  ball of radius $n$.
  (These and other pathological examples were pointed out to the
  originally over-optimistic author by C.\ Ward Henson.)
\end{rmk}

\begin{dfn}
  Say that a family of conditions $\Sigma = \{\varphi_i\leq r_i\colon i \in \lambda\}$ is
  \emph{approximately finitely satisfiable} if for every finite
  $w \subseteq \lambda$ and $\varepsilon > 0$, the family
  $\Sigma_0 = \{\varphi_i\leq r_i+\varepsilon\colon i \in w\}$ is satisfiable.
\end{dfn}

\begin{cor}
  If a set of sentential conditions (i.e., conditions without free
  variables) is approximately finitely satisfied in a family of
  structures, then it is satisfied
  in some ultra-product of these structures.
\end{cor}
\begin{proof}
  Standard.
\end{proof}

\begin{cor}[Bounded compactness for unbounded continuous logic]
  \label{cor:UnbddCpt}
  Let $\cL$ be an unbounded signature, $r \in \bR^+$,
  and let $\Sigma$ be a family of conditions in the free variables
  $x_{<n}$.
  Then $\Sigma \cup \{\nu(x_i) \leq r\colon i < n\}$ is satisfiable of and only if it is
  approximately finitely satisfiable.
\end{cor}

As usual, a theory is a set of sentential conditions.
The complete theory of a structure $M$, elementary equivalence and
elementary embeddings are defined as usual.
\begin{cor}
  Two structures $M$ and $N$ are elementarily equivalent if and only
  if $M$ embeds elementarily into an ultra-power of $N$.
\end{cor}
\begin{proof}
  One direction is clear.
  For the other we observe that if $M$ and $N$ are elementarily
  equivalent, then the elementary diagram of $M$ is approximately
  finitely satisfiable in $N$.
  Indeed, let $\bar a \in M$ and say that $\varphi(\bar a) = 0$.
  Let also $\varepsilon > 0$ and $r = \nu(\bar a)$.
  Then
  $N \models \inf_{\bar x}^{r,r+\varepsilon} \varphi(\bar x) = 0$,
  so there are $\bar b \in N$
  such that $\nu(\bar b) < r+\varepsilon$ and
  $\varphi(\bar b) < \varepsilon$.
\end{proof}

We could prove an analogue of the Shelah-Keisler theorem that if $N$
and $M$ are elementarily equivalent then they have isomorphic
ultra-powers.
We give a more elementary proof of a lesser result, which will
suffice just as well later on.
\begin{lem}
  \label{lem:CharElem}
  \begin{enumerate}
  \item 
    Two models $M$ and $N$ are elementarily equivalent if and only if
    there are sequences $M = M_0 \preceq M_1 \preceq \ldots$ and $N = N_0 \preceq N_1 \preceq \ldots$
    where each $M_{n+1}$ ($N_{n+1}$) is an ultra-power of $M_n$ ($N_n$)
    and $\bigcup_{n\in\bN} M_n \simeq \bigcup_{n\in\bN} N_n$
    (so their completions are isomorphic as well).
  \item A class of structures $\cK$ is elementary if and only if it is
    closed under elementary equivalence and ultra-products.
  \end{enumerate}
\end{lem}
\begin{proof}
  For the first item, right to left
  by the elementary chain lemma, which is proved as
  usual.
  For left to right, assume that $M \equiv N$.
  Then there is an ultra-power $N_1 = N^\sU$ and an elementary
  embedding $f_0\colon M \to N_1$.
  Then $(M,M) \equiv (N_1,f_0(M))$ (in a language with all elements of $M$
  named) so there exists an ultra-power $M_1 = M^{\sU'}$ and an
  elementary embedding
  $g_0\colon N_1 \to M_1$ such that $g_0 \circ f_0 = \id_M$.
  Proceed in this manner to obtain the sequences.

  The second item is standard.
\end{proof}

It is easily verified that any theory is logically equivalent to one
which only consists of conditions of the form $\varphi = 0$.
A \emph{universal} theory is one which only consists
of conditions of the form $\sup_{\bar x} \varphi(\bar x) = 0$ where $\varphi$ is
quantifier-free
(and syntactically bounded and eventually constant in each $x_i$).
Observe that:
\begin{itemize}
\item For any formula $\varphi$ we can express $\forall\bar x\,
  \varphi(\bar x) = 0$
  by the universal axiom scheme
  $\sup_{\bar x}^{n,n+1} \varphi(\bar x) = 0$.
\item If $t$ and $s$ are terms we can express $\forall\bar x\, t = s$
  by $\forall\bar x\, d(t,s) = 0$.
\item If $\varphi$ and $\psi$ are formulae we can express
  $\forall\bar x\, \varphi \geq \psi$ by
  $\forall\bar x\, \psi \dotminus \varphi = 0$.
\end{itemize}

\begin{exm}
  \label{exm:BanachTh}
  We can continue \fref{exm:BanachSig} and give the (universal)
  theory of the class of Banach spaces:
  \begin{align*}
    &
    \langle
    \textit{Universal equational axioms of a vector space.}
    \rangle
    \\ &
    \forall x\, s\|x\| \leq \|m_r(x)\| \leq s'\|x\|
    && s,s' \text{ dyadic}, s \leq |r| \leq s'
    \\ &
    \forall xy \, \|x+y\| \leq \|x\| + \|y\|
    \\ &
    \forall xy \, d(x,y) = \|x + m_{-1}(y)\|.
  \end{align*}
\end{exm}

More generally, it will be convenient to write
\begin{gather*}
  \bigl( {\sup_x}^r \, {\inf_y}^s \, \ldots \, \varphi \bigr) = 0,
  \qquad \text{or even} \qquad
  \forall^{<r}x \, \exists^{\leq s}y \, \ldots \bigl( \varphi = 0 \bigr)
\end{gather*}
for the axiom scheme
\begin{gather*}
  {\sup_x}^{r-\varepsilon,r} {\inf_y}^{s,s+\varepsilon} \ldots
  \varphi = 0,
  \qquad \varepsilon > 0.
\end{gather*}
Notice that in the $\forall\exists$ notation, the universal quantifier
holds literally, while the existential quantifiers holds in an
approximate sense, with respect to the quantification radius as well
as with respect to the value of $\varphi$
(which may both be slightly bigger than $s$ or $0$, respectively.)

\begin{exm}[Measure algebras]
  \label{exm:MeasureAlgebras}
  Let $\cL = \{0,\vee,\wedge,\setminus\}$, where $0$ is a constant symbol, $\vee,\wedge,\setminus$ are
  binary function symbols.
  We use $\mu(x)$ as shorthand for $d(x,0)$, and take it to be the
  gauge.
  All the continuity moduli are standard.

  The universal theory of measure algebras (which are the topic of
  \cite{Fremlin:MeasureTheoryVol3}) consists of:
  \begin{align*}
    & \langle\textit{Universal equational axioms of relatively
      complemented distributive lattices}\rangle, \\
    & \forall xy \, \mu(x) + \mu(y) = \mu(x\wedge y) + \mu(x\vee y), \\
    & \mu(0) = 0, \\
    & \forall xy\, d(x,y) = \mu(x\setminus y) + \mu(y\setminus x).
  \end{align*}

  We can further say that a measure algebra is \emph{atomless} by the
  axiom scheme:
  \begin{align*}
    & \forall^{<n}x \, \exists^{\leq n}y \,
    |\mu(x \wedge y) - \mu(x)/2| = 0,
    && n \in \bN.
  \end{align*}
\end{exm}

\begin{exm}[Replacing a function with its graph]
  \label{exm:FunctionGraph}
  Let $\cL$ be an unbounded signature, $f \in \cL$ an $n$-ary function
  symbol.
  We define its \emph{graph} to be the $(n+1)$-ary predicate
  $G_f(\bar x,y) = d(f(\bar x),y)$.
  Since it is defined by a formula it respects a continuity
  modulus under $\nu$ uniformly in all $\cL$-structures, and we may
  add it to the language.
  The axiom scheme
  $\forall\bar x y\, G_f(\bar x,y) = d(f(\bar x),y)$ is universal.

  We may further drop $f$ from the language.
  Indeed, we observe that a predicate $G_f$ is the graph of a function
  $f$ with continuity modulus $\delta_f$ if and only if the following
  theory holds.
  The second axiom ensures that in the third axiom there actually
  exists a unique $y = f(\bar x)$ such that $G_f(\bar x,y) = 0$.
  Then the first two axioms imply that $G_f$ is the graph of $f$,
  and the two last axioms together ensure that
  $f$ respects $\delta_f$ under $\nu$.
  \begin{align*}
    & \forall \bar x,y,z  \,
    && G_f(\bar x,y) \leq G_f(\bar x,z) + d(y,z)
    \\
    & \forall \bar x,y,z  \,
    && d(y,z)  \leq G_f(\bar x,y) + G_f(\bar x,z)
    \\
    &
    \forall^{<\varepsilon^{-1}}\bar x
    \,
    \exists^{\leq\delta_f(\varepsilon)^{-1}}y
    \,
    &&G_f(\bar x,y) = 0
    && \varepsilon > 0
    \\
    &
    \forall^{<\varepsilon^{-1}}\bar x\bar y
    \,
    \forall^{<\delta_f(\varepsilon)^{-1} +1}z
    \,
    &&
    \bigl( \delta_f(\varepsilon) \dotminus d(\bar x,\bar y) \bigr)
    \wedge
    \bigl( G_f(\bar x,z) \dotminus G_f(\bar y,z) \dotminus \varepsilon \bigr)
    =0
    && \varepsilon > 0
  \end{align*}
\end{exm}

Types and type spaces are defined more or less as usual:
\begin{dfn}
  Fix an unbounded signature $\cL$.
  \begin{enumerate}
  \item Given an $n$-tuple $\bar a$, we define its type
    $p(\bar x) = \tp(\bar a)$ as usual
    as the set of all $\cL$-conditions in the variables $x_{<n}$
    satisfied by $\bar a$.
    The type $p(\bar x)$ determines the value of $\varphi(\bar a)$ for
    every formula $\varphi$, and we may write
    $\varphi^p = \varphi(\bar x)^{p(\bar x)} = \varphi(\bar a)$.
  \item A complete $n$-type (in $\cL$) is the type of some $n$-tuple.
    By \fref{cor:UnbddCpt}, this is the same as a maximal 
    finitely consistent set of conditions $p(x_{<n})$
    such that for some $r \geq 0$ we have $\nu(x_i) \leq r \in p$ for all
    $i < n$.
  \item The set of all $n$-types is denoted $\tS_n$.
    The set of all $n$-types containing a theory $T$ (equivalently:
    realised in models of $T$) is denoted $\tS_n(T)$.
  \item For every condition $s$ in the free variables $x_{<n}$,
    $[s]^{\tS_n(T)}$ (or just $[s]$, if the ambient type space is
    clear from the context) denotes the set of types
    $\{p \in \tS_n(T)\colon s \in p\}$.
  \item The family of all sets of the form $[s]^{\tS_n(T)}$ forms a
    base of closed sets for the \emph{logic topology} on $\tS_n(T)$.
    It is easily verified to be Hausdorff.
  \end{enumerate}
\end{dfn}

For each $n \in \bN$, we can define $\nu\colon \tS_n(T) \to \bR$ by
$\nu(p) = \bigvee_{i<n} \nu(x_i)^p$.
With this definition, $(\tS_n(T),d,\nu)$ is a gauged space.
Applying previous definitions we have:
\begin{gather*}
  \tS^{\nu\leq r}_n(T)
  = \bigcap_{i<n} \bigl[ \nu(x_i) \leq r \bigr]
  = \left[ \Bigl( \bigvee_{i<n} \nu(x_i) \Bigr) \leq r \right].
\end{gather*}
By \fref{cor:UnbddCpt}, $\tS^{\nu\leq r}_n(T)$ is compact.
If $\tS_n(T) = \tS^{\nu\leq r}_n(T)$ for some $r$, then $\tS_n(T)$ is compact.
Conversely, if $\tS_n(T)$ is compact for $n \geq 1$, then $\nu$ is
necessarily bounded on models of $T$, so there is some $r$ such that
$T \vdash \sup_x \nu(x)\wedge(r+1) \leq r$ and
$\tS_m(T) = \tS_m^{\nu\leq r}(T)$ for all $m \in \bN$.
In this case all the other symbols are also bounded in models
of $T$, so up to re-scaling everything
into $[0,1]$ we are in the case of
standard continuous first order logic.

In the non compact case we still have
$\tS_n(T) = \bigcup_r \tS^{\nu\leq r}_n(T)$.
Thus each $p \in \tS_n(T)$ there is $r$ such that
$p \in \tS^{\nu\leq r}_n(T)$, and $\tS_n^{\nu\leq r+1}(T)$
is a compact neighbourhood of $p$
(since it contains the open set $[(\bigvee\nu(x_i)) < r+1]$).
Therefore $\tS_n(T)$ is locally compact.

\section{On the relation with Henson's positive bounded logic}
\label{sec:Henson}

We sketch out here how unbounded continuous logic generalises,
in an appropriate sense, Henson's logic of
approximate satisfaction of positive bounded formulae in Banach space
structures.
For this purpose we assume familiarity with the syntax and
semantics of Henson's logic
(see for example \cite{Henson-Iovino:Ultraproducts}).

The classical presentation of Henson's logic involves a purely
functional signature $\cL_H$ with a distinguished sort for $\bR$.
There is no harm in assuming that the distinguished sort
only appears as the target sort of some function symbols
(otherwise we can add a second copy and a single function symbol for
the identity mapping into the copy, and treat the copy as the
distinguished sort).
Also, there is no harm in replacing $\bR$ with $\bR^+$.

We can therefore define an unbounded continuous signature $\cL$
by dropping the distinguished sort and replacing all function
symbols into it with $\bR^+$-valued predicate symbols.
As every sort is assumed to be normed, we identify $\nu$ with
$\|{\cdot}\|$.
While a signature in Henson's logic does not specify continuity
moduli, in every class under consideration each symbol satisfies some
continuity modulus uniformly under $\|{\cdot}\|$ which we may use
(or else the logic would fail to describe the class).
It is a known fact that there exists a (universal)
$\cL_H$-theory, call it $T_0$, whose models are precisely the
structures respecting these continuity moduli under $\|{\cdot}\|$.

From now on by ``structure'' we mean a model of $T_0$,
or equivalently a $\cL$-structure (as these can be identified).
The ambiguity concerning whether a structure is a Henson or unbounded
continuous structure is further justified by the fact that
the definitions of isomorphism and ultra-products in either logic
coincide.
As we can moreover prove \fref{lem:CharElem} for Henson's logic
just as well, we conclude:
\begin{thm}
  \label{thm:EqUnbddHenson}
  A class of structures $\cK$ is elementary in Henson's logic if and
  only if it is elementary in unbounded continuous logic.
\end{thm}

Recall:
\begin{fct}
  \label{fct:UnbddClsd}
  Let $X = \bigcup_{n\in\bN} X_n$ be a topological space where each
  $X_n$ is closed and $X_{n+1}$ is a neighbourhood of $X_n$.
  Then a subset $F \subseteq X$ is closed if and only if $F\cap X_n$ is for all
  $n$.
\end{fct}

An $n$-type is the same thing as a complete theory with
$n$ new constant symbols
(more precisely, a type $p$ with $\nu(p) \leq r$ corresponds
to a complete theory with new constants symbols with continuity
moduli $\delta_a \leq \inv r$).
\begin{cor}
  Two $n$-tuples in a structure have the same type in one logic if and
  only if they have the same type in the other, and this
  identification induces a homeomorphism
  $\tS_n^{\cL_H}(T_0) \simeq \tS_n^\cL$.
\end{cor}
\begin{proof}
  The first statement is by \fref{thm:EqUnbddHenson}.
  Also, a set $X \subseteq \tS_n^{\|\cdot\|\leq r}$ is closed if and only if
  the class $\{(M,\bar a)\colon \tp(\bar a) \in X\}$ is elementary: the bounds
  on the norm are needed since we need to impose bounds on the norms
  of constant symbols.
  It follows from \fref{thm:EqUnbddHenson} that
  the bijection $\tS_n^{\cL_H}(T_0) \simeq \tS_n^\cL$ is a homeomorphism
  when restricted to $\tS_n^{\|\cdot\|\leq n}$.
  Now use \fref{fct:UnbddClsd} and the fact that $\tS_n^{\|\cdot\|\leq r}$ is
  compact and $\tS_n^{\|\cdot\|<r}$ is open in both topologies
  to conclude that this is a global homeomorphism.
\end{proof}

This can be restated as:
\begin{cor}
  For every set $\Sigma(\bar x)$ of $\cL_H$-formulae there
  exists a set $\Gamma(\bar x)$ of $\cL$-conditions, and for
  every set $\Gamma(\bar x)$ of $\cL$-conditions
  there exists a set $\Sigma(\bar x)$ of $\cL_H$-formulae, such
  that for every structure $M$ and $\bar a \in M$:
  \begin{gather*}
    M \models_A \Sigma(\bar a)
    \quad \Longleftrightarrow \quad
    M \models \Gamma(\bar a).
  \end{gather*}
\end{cor}

\begin{rmk}
  In Henson's logic, the bounded quantifier $\forall^{\leq r}x$
  ($\exists^{\leq r}x$)
  mean ``for all (there exists) $x$ such that $\|x\|\leq r$''.
  Thus Henson's logic coincides with unbounded continuous logic of
  normed structures where $\nu = \|{\cdot}\|$.
  One may generalise Henson's logic to allow an arbitrary $\nu$ and
  obtain full equivalence of the two logics.
\end{rmk}

For the benefit of the reader who finds this proof a little too
obscure, let us give one direction explicitly.
We know that every formula in Henson's logic is equivalent to one in
prenex form
\begin{align*}
  & \forall^{\leq r_0} x_0 \exists^{\leq r_1}x_1 \ldots \,
  \varphi(\bar x,\bar y),
\end{align*}
where $\varphi$ is a positive Boolean combination of atomic formulae
of the form $t_i(\bar x,\bar y) \geq r_i$ or $t_i \leq r_i$.
Every term $t_i$ can be identified with an atomic $\cL$-formula,
and replacing $t_i$ with $t_i \dotminus r_i$ or with $r_i \dotminus t_i$,
we may assume all these atomic formulae are of the form
$t_i \leq 0$.
Since $(t_i\leq0)\wedge(t_j\leq0) \Longleftrightarrow (t_i\vee t_j)\leq0$ and
$(t_i\leq0)\vee(t_j\leq0) \Longleftrightarrow (t_i\wedge t_j)\leq0$, we can find a single $t$ such that
$\varphi(\bar x,\bar y)$ is equivalent to $t \leq 0$.
We thus reduced to:
\begin{align*}
  & \forall^{\leq r_0} x_0 \, \exists^{\leq r_1}x_1 \ldots \,
  \bigl( t(\bar x,\bar y) \leq 0 \bigr).
\end{align*}
We can view $t$ as a quantifier-free $\cL$-formula, in which case the
above holds approximately if and only if the following holds
(with the notation preceding \fref{exm:MeasureAlgebras}):
\begin{align*}
  & \forall^{< r_0} x_0 \, \exists^{\leq r_1}x_1 \ldots \,
  t(\bar x,\bar y) = 0.
\end{align*}
Thus the approximate satisfaction of a $\cL_H$-formula, and therefore
of a partial type, are equivalent to the satisfaction of
a partial type in $\cL$.

\section{Emboundment}
\label{sec:Emboundment}

As we mentioned earlier, the multi-sorted approach to unbounded
structures allows us to reduce many issues concerning unbounded
structures to their well-established analogues in bounded continuous
logic, but this does not work well for
perturbations when we wish to perturb $\nu$ itself.
In addition, if the bounded balls are not definable in the
unbounded structure then their introduction as sorts adds
unexpected structure -- this may happen, for example, when considering
a field equipped with a valuation in $\bR$ as an unbounded metric
structure.

We could of course generalise everything we did to the unbounded case,
but that would be extremely tedious to author and reader alike.
Instead, we seek a universal reduction of unbounded logic to the
more familiar (and easier to manipulate) bounded one.
This reduction goes through a construction which we call
\emph{emboundment}.
Thus, for example, a bounded set $X \subseteq M^n$ in an unbounded
structure is said to be \emph{definable} (a term we knowingly used
above without a definition) if it is definable in the embounded
structure $M^\infty$.
An easy verification yields that this is equivalent to the predicate
$d(\bar x,X)$ being definable in $M$, i.e., a uniform limit of
formulae on every bounded set.
(See \cite{BenYaacov:DefinabilityOfGroups} for definable sets in
bounded structures.)

One naïve approach would be to choose a continuous function mapping
$\bR^+$ into $[0,1]$, say $\theta(x) = \frac{x}{x+1}$, and apply it to all
the predicate symbols: for every $\cL$-structure $M$ we define
$M^\theta$ as having the same underlying set, and for every predicate
symbol $P$ we define $P^{M^\theta}(\bar a) = \theta(P^M(\bar a))$.
It can be verified that
$\theta(x+y) \leq \theta(x)+\theta(y)$
for all $x,y \geq 0$
(this is true when $x = 0$, and the partial derivative with
respect to $x$ of the left hand side is smaller).
It follows that $d^{M^\theta}$ is a metric:
\begin{align*}
  d^{M^\theta}(a,b) &
  = \theta(d^M(a,b)) \leq \theta(d^M(a,c) + d^M(c,b))
  \\ &
  \leq \theta(d^M(a,c)) + \theta(d^M(c,b))
  = d^{M^\theta}(a,b) + d^{M^\theta}(a,b).
\end{align*}
Of course $d^{M^\theta}$ needs not be a complete metric, so we obtain new
elements when passing to the completion.
Similarly, if $T^\theta = \Th\{M^\theta\colon M \models T\}$,
then we have a natural embdding
of $\tS_n(T)$ in $\tS_n(T^\theta)$, and it can be verified that the latter
is the Stone-Čech compactification of the former.
This is essentially the same thing as allowing $\infty$ as a legitimate
truth value (as $\theta$ extends to a homeomorphism $[0,\infty] \to [0,1]$).
As usual with the Stone-Čech compactification, this adds too many new
types to be manageable.
In short, this naïve construction does yield bounded structures
but it is not at all clear that the structures (or theories) thus
obtained are meaningful.
For example, even the following is not clear (to the author), and one
would expect it to be false:
\begin{qst}
  Is every model of $T^\theta$ of the form $M^\theta$, where $M \models T$?
\end{qst}

For a better approach, we take a second look on the construction of
unbounded logic and its semantics, as well as on the construction of
unbounded ultra-products.
Throughout these constructions appeared a formal infinity element
$\infty$, which, while not a member of the structures, was treated for
many intents and purposes as if it were.
Indeed, the quantifier semantics included $\infty$ in the set over
which quantification takes place, and the ultra-product construction
could be restated informally as ``add $\infty$, take a usual
ultra-product, then take $\infty$ out''.
In particular, unbounded structures may be formally empty since, from
a practical point of view, they still always contain the ideal point
at infinity.

With this motivation in mind, we seek to equip each unbounded
structure $M$ with a new metric, denoted $d^{M^\infty}$
such that
every sequence $(a_n)$ in $M$ which goes to infinity in the sense that
$\nu(a_n) \to \infty$, is Cauchy in $d^{M^\infty}$, converging to a
new element representing the formal infinity.
Such a metric is naturally bounded.
Moreover, every predicate on $M$ which is uniformly continuous under
$\nu$ can be modified to yield a bounded predicate which is in
uniformly continuous in the usual sense with respect to
$d^{M^\infty}$.
On the other hand, this does not work well for function symbols
(for example, we cannot give a sense to $\infty + \infty$ in the
emboundment of a Banach space).
We shall therefore replace every function symbol in the language with
its graph $G_f(\bar x,y) = d(f(\bar x),y)$ as in
\fref{exm:FunctionGraph}, and assume that the signature $\cL$ is
purely relational.
We then define
\begin{gather*}
  \cL^\infty = \cL \cup \{\infty\}
\end{gather*}
where $\infty$ is a new constant symbol.
We may consider $\cL$ to consist, as a set, of its non distinguished
symbols alone, in which case $\nu$ gets dropped
(or more precisely, both $d$ and $\nu$ are dropped, and then
$\cL^\infty$ is equipped with its own distinguished distance symbol
$d$).
Whether or not $\nu$ is kept will be of no essential difference to the
construction.
We do not specify at this point the uniform continuity moduli, but
we shall show below that such moduli can be chosen that do fit our
purpose.

For every $\cL$-structure $M$ we define an $\cL^\infty$-structure
$M^\infty$.
Its domain is the set $M \cup \{\infty\}$.
For elements coming from $M$ we interpret the symbols as follows
(we recall that $d(a,\infty) = \nu(\infty) = \infty$,
$\theta(\infty) = 1$, and $\nu(\bar x) = \bigvee \nu(x_i)$):
\begin{align*}
  &
  d^{M^\infty}(a,b)
  = \frac{\theta\circ d^M(a,b)}{1+\nu^M(a) \wedge \nu^M(b)},
  &&
  P^{M^\infty}(\bar a)
  = \frac{\theta \circ P^M(\bar a)}{1+\nu^M(\bar a)},
  &&
  (P \neq d).
  \intertext{So in particular:}
  &
  d^{M^\infty}(a,\infty)
  = \frac{1}{1+\nu^M(a)},
  &&
  P^{M^\infty}(\ldots,\infty,\ldots) = 0,
  &&
  (P \neq d).
\end{align*}
Notice that if we interpreted $d^{M^\infty}$ as with other symbols we
would have $d^{M^\infty}(a,\infty) = 0$ for all $a$,
and thus not obtain a metric.
Conversely, we can reconstruct $M$ from $M^\infty$,
first recovering $\nu^M$ from $d^{M^\infty}(x,\infty)$ and
then recovering $d^M$ and $P^M$ from $d^{M^\infty}$
and $P^{M^\infty}$, respectively,
using the fact that
$\theta^{-1}(y) = \frac{y}{1-y}$.

Let us show that $d^{M^\infty}$ is a metric.
The only non trivial property to verify is the triangle inequality,
namely
\begin{gather*}
  \frac{\theta\circ d^M(a,c)}{1+\nu^M(a) \wedge \nu^M(c)}
  \leq
  \frac{\theta\circ d^M(a,b)}{1+\nu^M(a) \wedge \nu^M(b)}
  +
  \frac{\theta\circ d^M(b,c)}{1+\nu^M(b) \wedge \nu^M(c)}.
\end{gather*}
If $b$ has the smallest gauge among the three then this follows from
the fact that $\theta \circ d^M$ is a metric, which we verified
earlier.
Otherwise we may assume without loss of generality that $a$ has the
smallest gauge, say $r$.
Let $t = d^M(a,b)$, $s = d^M(b,c)$.
Then $\nu^M(b) \leq r+t$ and $d^M(a,c) \leq t+s$,
and it is enough to verify that
\begin{gather*}
  \frac{\theta(t+s)}{1+r}
  \leq
  \frac{\theta(t)}{1+r}
  +
  \frac{\theta(s)}{1+r+t}.
\end{gather*}
Moving the second term to the left and developing we obtain:
\begin{align*}
  \frac{\theta(t+s)}{1+r}
  -
  \frac{\theta(t)}{1+r}
  &
  =
  \frac{s}{(1+r)(1+t)(1+t+s)}
  \\ &
  =
  \frac{s}{(1+r)(1+t+s) + (1+r)t(1+t+s)}
  \\ &
  \leq
  \frac{s}{(1+r)(1+s) + t(1+s)}
  =
  \frac{s}{(1+r+t)(1+s)}
  \\ &
  =
  \frac{\theta(s)}{1+r+t},
\end{align*}
as desired.
Once we know that $d^\infty$ is a metric it is clear that
$a_n \to \infty$ in $d^{M^\infty}$ if and only if
$\nu^M(a_n) \to \infty$.

\begin{exm}
  Let $M$ be a bounded structure,
  and turn it into an unbounded structure $M'$ as in
  \fref{exm:BoundedStructreUnbounded}.
  Then $M \cong (M')^\infty \setminus \{\infty\}$,
  so all we did was add a single isolated point with distance $1$ to
  the original structure.
\end{exm}

\begin{lem}
  \label{lem:EmboundUC}
  The gauged space $(M,d^M,\nu^M)$ and the bounded metric space
  $(M,d^{M^\infty})$ are related as follows:
  \begin{enumerate}
  \item \label{item:EmboundUCMetrics}
    We have $d^M \geq d^{M^\infty}$ on all of $M$,
    and the two metrics are uniformly equivalent on every bounded
    subset of $M$ (bounded in the sense of $M$).
  \item \label{item:EmboundUCBalls}
    For every $r'>r$ the $\nu$-ball
    $M^{\nu< r'}$ contains a uniform $d^{M^\infty}$-neighbourhood of
    $M^{\nu\leq r}$ (of radius $\frac{\theta(r'-r)}{1+r}$).
  \end{enumerate}
\end{lem}
\begin{proof}
  The inequality $d^M \geq d^{M^\infty}$ is immediate.
  Let us fix $r \geq 0$ and let $a \in M^{\nu\leq r}$, $b \in M$.
  Then by definition
  $d^{M^\infty}(a,b) \geq \frac{\theta \circ d^M(a,b)}{1+r}$.
  Thus, for all $\varepsilon > 0$
  \begin{gather*}
    d^{M^\infty}(a,b) < \frac{\theta(\varepsilon)}{1+r}
    \quad \Longrightarrow \quad
    d^M(a,b) < \varepsilon,
  \end{gather*}
  concluding the proof of the first item.
  This also proves the third item, since

  \begin{gather*}
    B_{d^{M^\infty}}\left(
      M^{\nu\leq r},\frac{\theta(r'-r)}{1+r}
    \right)
    \subseteq
    B_d(M^{\nu\leq r},r'-r)
    \subseteq M^{\nu<r'}.
    \qedhere
  \end{gather*}
  
\end{proof}

\begin{prp}
  For every $\cL$-structure $M$, $M^\infty$ as defined above is an
  $\cL^\infty$-structure, called the \emph{emboundment} of $M$.
  That is to say that $M^\infty$ is complete, and
  that we can complete the definition of $\cL^\infty$ choosing uniform
  continuity moduli for its symbols which are satisfied in every
  $M^\infty$.
\end{prp}
\begin{proof}
  For completeness, let $(a_n)_{n\in\bN}$ be a Cauchy sequence in
  $d^{M^\infty}$.
  If $\nu(a_n) \to \infty$ (where again, $\nu(\infty) = \infty$)
  then $a_n \to \infty$ in $d^\infty$.
  Otherwise, there is $r$ such that $a_n \in M^{\nu< r}$ infinitely
  often.
  Passing to a sub-sequence, we may assume that the entire sequence fits
  inside $M^{\nu<r}$.
  By \fref{lem:EmboundUC}\fref{item:EmboundUCMetrics} the sequence is
  Cauchy in $d^M$ and therefore admits a limit in $M$, which is
  necessarily also its limit in $d^{M^\infty}$.

  For uniform continuity, let
  $P \in \cL$ be an $(n+1)$-ary predicate symbol.
  Let $\varepsilon > 0$ be given, and we wish to find $\delta > 0$ such that
  for all $a,b \in M^\infty$
  \begin{gather*}
    d^{M^\infty}(a,b) \leq \delta
    \quad \Longrightarrow \quad
    \sup_{\bar x} |\tilde P(a,\bar x) - \tilde P(b,\bar x)|^{M^\infty}
    \leq \varepsilon.
  \end{gather*}
  First, if $\nu(a),\nu(b) \geq \inv\varepsilon-1$
  (where $\nu(\infty) = \infty$)
  then the above is satisfied regardless of
  $d^{M^\infty}(a,b)$.
  Otherwise, without loss of generality we have
  $\nu(a) < \inv\varepsilon-1$, and if
  $d^{M^\infty}(a,b) < \frac{\varepsilon}{2}$
  then $\nu(b) < \inv\varepsilon$
  by \fref{lem:EmboundUC}\fref{item:EmboundUCBalls}.
  Since $P$ and $\nu$ are uniformly continuous
  with respect to $d^M$ on $M^{\nu<\frac{1}{\varepsilon}}$,
  so is $\tilde P$.
  By \fref{lem:EmboundUC}\fref{item:EmboundUCMetrics},
  $\tilde P$ is uniformly continuous with respect to
  $d^{M^\infty}$ on $M^{\nu<\frac{1}{\varepsilon}}$,
  whence the existence of $\delta$ as desired.
\end{proof}

It is straightforward to verify that the emboundment
construction commutes with the ultra-product construction, since
everything is continuous:
$$\left( \prod M_i/\sU \right)^\infty  = \prod M_i^\infty/\sU.$$
In particular, all the tuples $(a_i)$ such that
$\lim_\sU \nu^{M_i}(a_i) = \infty$, which were dropped during the
construction of
$\prod M_i/\sU$, satisfy $[a_i] = [\infty^{M_i}] = \infty$ in $\prod
M_i^\infty/\sU$.

Similarly, emboundment commutes with unions of increasing chains,
and by \fref{lem:CharElem} we have
$M \equiv N \Longleftrightarrow M^\infty \equiv N^\infty$
for any two $\cL$-structures $M$ and $N$.
If $N \subseteq M$ then working with $\cL(N)$ we get
$N \preceq M \Longleftrightarrow N^\infty \preceq M^\infty$.
Similarly, if
$N' \preceq M^\infty$ where $N'$ is an $\cL^\infty$-structure
then we can recover an $\cL$-structure
on $N = N' \setminus \{\infty\} \subseteq M$,
so $N' = N^\infty$ and $N \preceq M$.

\begin{prp}
  \label{prp:SPCElem}
  Let $\cK$ be a class of $\cL$-structures, and let
  $$\cK^\infty = \{M^\infty\colon M \in \cK\}.$$
  Then $\cK$ is elementary if and only if $\cK^\infty$ is.
\end{prp}
\begin{proof}
  Assume $\cK$ is elementary.
  Then, by the arguments above, $\cK^\infty$ is closed under
  ultra-products, isomorphism and elementary substructures.
  It is therefore elementary.
  Similarly for the converse.
\end{proof}

By \fref{prp:SPCElem} we may replace every
$\cL$-theory $T$ (in unbounded logic) with its emboundment
$T^\infty = \Th_{\cL^\infty}(\Mod(T)^\infty)$, which is a theory in
standard bounded logic.
By naming constants we further see
that $\cL$-types of tuples in $M$
are in bijection with $\cL^\infty$-types of tuples in
$M^\infty \setminus \{\infty\}$
(i.e., in $M$ again, but this time viewed as a subset of a
$\cL^\infty$-structure).

Given a tuple $\bar a \in M^\infty$, let
$w = w(\bar a) = \{i < n\colon a_i \neq \infty\}$.
Then we may identify $\tp^{M^\infty}(\bar a)$
with the pair $(w,\tp^M(a_{\in w}))$.
We can therefore express the set of types $\tS_n(T^\infty)$ as
$\bigcup_{w\subseteq n} \{w\}\times\tS_{|w|}(T)$.
For $w \subseteq n$, $r \in \bR^+$ and $\varphi(x_{\in w}) \in \cL$, define:
\begin{align*}
  V_{n, w,r,\varphi} = \left\{ (v,q(x_{\in v})) \in \tS_n(T^\infty)\colon
  \begin{aligned}
    & w \subseteq v \subseteq n,\; \varphi^q < 1, \\
    & {\bigwedge}_{i \in v \setminus w}\nu(x_i)^q > r
  \end{aligned}
  \right\}.
\end{align*}
Given a type $(w,p) \in \tS_n(T^\infty)$, one can verify that the family of
all sets of the form $V_{n,w,r,\varphi}$ where $\varphi^p = 0$ forms a base of
neighbourhoods for $(w,p)$.
In particular, the natural inclusion $\tS_n(T) \hookrightarrow \tS_n(T^\infty)$,
consisting of sending $p \mapsto (n,p)$, is an open topological embedding.
In case $T$ is complete (so $|\tS_0(T)| = 1$), this embedding for
$n = 1$ is a single point compactification of $\tS_1(T)$
obtained by adding the type at infinity.
We may therefore also refer to $T^\infty$ as the
\emph{compactification} of $T$.

Once we understand types we know what saturation means.
Among other things we have:
\begin{lem}
  \label{lem:UnbddOSat}
  An $\cL$-structure $M$ is approximately $\aleph_0$-saturated if and only if
  $M^\infty$ is.
\end{lem}
\begin{proof}
  Follows from the facts that there is a unique point at infinity,
  which belongs to $M^\infty$, and that in
  the neighbourhood of every other point $d^M$ and
  $d^{M^\infty}$ are equivalent.
\end{proof}

Finally, we point out that
the theory $T$ is bounded to begin with if and only if
the point at infinity in models of $T^\infty$ is
isolated, in analogy with what happens when
one attempts to add a point at infinity to a space which is already
compact.

\section{Perturbations of unbounded structures}
\label{sec:Perturbation}

We now adapt the framework of perturbation of bounded metric
structures to unbounded structures,
essentially by reducing the unbounded case to the bounded one through
emboundment.
For this purpose we assume close familiarity with the original
development in \cite{BenYaacov:Perturbations}.
We fix an unbounded theory $T$ and its emboundment $T^\infty$.

\begin{dfn}
  A \emph{perturbation pre-radius} for
  $T$ is defined as for a bounded theory, i.e., as a
  family $\rho = \{\rho_n \subseteq \tS_n(T)\colon n \in \bN\}$ containing the diagonals.
  We define $X^\rho$, $\Pert_\rho(M,N)$, $\BiPert_\rho(M,N)$,
  $\langle\rho\rangle$, $\llbracket\rho\rrbracket$ as in \cite{BenYaacov:Perturbations}.
\end{dfn}

Let $\rho$ be a perturbation pre-radius for $T$.
We can always extend it to a perturbation radius $\rho^*$ for $T^\infty$ by:
\begin{align*}
  \rho^\infty_n & = \overline{\big\{
    ((w,p),(w,q)) \in \tS_n(T^\infty)\colon w \subseteq n, (p,q) \in \rho_{|w|}
    \big\}}.
\end{align*}
Clearly, this is a perturbation pre-radius for $T^\infty$.
Conversely, if $\rho'$ is a perturbation pre-radius for $T^\infty$ then its
restriction to $\tS(T)$, denoted $\rho'\rest_{\tS(T)}$, is a perturbation
pre-radius for $T$, and as the inclusion $\tS_n(T) \subseteq \tS_n(T^\infty)$ is
open we have the identity:
$$\rho^\infty\rest_{\tS(T)} = \rho.$$
Also, as every $f \in \Pert_\rho(M,N)$ extends to
$f\cup(\infty\mapsto\infty) \in \Pert_{\rho^\infty}(M^\infty,N^\infty)$, we also have $\langle\rho^\infty\rangle\rest_{\tS(T)} \geq \langle\rho\rangle$.

We define perturbation radii for $T$ directly by reduction to $T^\infty$:
\begin{dfn}
  \begin{enumerate}
  \item Let $\rho'$ be a perturbation pre-radius for $T^\infty$.
    We say that $\rho'$ \emph{separates infinity} if for all
    $f \in \Pert_{\rho'}(M^\infty,N^\infty)$ and $a \in M^\infty$:
    $$a = \infty \Longleftrightarrow f(a) = \infty.$$
  \item A perturbation pre-radius $\rho$ for $T$ is a
    \emph{perturbation radius} if $\rho^\infty$ is a perturbation radius for
    $T^\infty$ which separates infinity.
  \end{enumerate}
\end{dfn}

\begin{dfn}
  A \emph{perturbation pre-system}
  for $T$ is a decreasing family $\fp$ of 
  perturbation pre-radii satisfying downward continuity, symmetry,
  triangle inequality and strictness as in
  \cite[Definition~1.23]{BenYaacov:Perturbations}.
  It is a \emph{perturbation system} if $\fp(\varepsilon)$ is a
  perturbation radius for all $\varepsilon$, i.e., if $\fp^\infty$
  is a perturbation system separating infinity for $T^\infty$.
\end{dfn}

We turn to characterise perturbation radii as in
\cite{BenYaacov:Perturbations}, and establish more precisely
the relation between perturbations of $T$ and of $T^\infty$.
\begin{dfn}
  Let $\rho$ a perturbation pre-radius for $T$.
  \begin{enumerate}
  \item We say that $\rho$ \emph{respects infinity} if for all
    $r \in \bR^+$ there exists $r' \in \bR^+$ such that
    \begin{gather*}
      [\nu(x)\geq r']^\rho \subseteq [\nu(x)\geq r]
      \qquad \text{ and } \qquad
      [\nu(x)\leq r]^\rho \subseteq [\nu(x)\leq r'].
    \end{gather*}
  \item We define when $\rho$ \emph{respects equality},
    \emph{respects $\exists$},
    or is \emph{permutation-invariant} as in the bounded
    case.
  \end{enumerate}
\end{dfn}

\begin{prp}
  \label{prp:UnbddPertRadFunct}
  Let $\rho$ be a perturbation pre-radius for $T$.
  The the following are equivalent:
  \begin{enumerate}
  \item $\rho$ is a perturbation radius.
  \item $\rho$ respects infinity, and for every $n,m \in \bN$
    and mapping $\sigma\colon n \to m$, the induced mapping
    $\sigma^* \colon \tS_m(T) \to \tS_n(T)$ satisfies
    that for all $p \in \tS_m(T)$:
    $$\sigma^*(p^\rho) = \sigma^*(p)^\rho.$$
    (I.e., $\sigma^* \circ \rho_m = \rho_n \circ \sigma^*$ as multi-valued functions).
  \item $\rho$ respects $\infty$, $=$, $\exists$, and is permutation-invariant.
  \item $\rho^\infty$ separates $\infty$, respects $=$ and $\exists$ and is
    permutation-invariant.
  \end{enumerate}
\end{prp}
\begin{proof}
  \begin{cycprf}
  \item[\impnext]
    Assume $\rho$ is a perturbation radius, so $\rho^\infty$ is a perturbation
    radius respecting infinity.
    If $\rho$ does not respect infinity, then by definition of $\rho^\infty$ we
    have in $\rho^\infty_1$ a pair $(p,q)$ where $p$ is the type of a
    finite elements and $q = \tp(\infty)$ or vice versa, contradicting the
    assumption on $\rho^\infty$.

    Since $\rho^\infty$ is a perturbation radius,
    for all $\sigma\colon n\to m$ we have in $\tS(T^\infty)$:
    $\sigma^* \circ \rho_m^\infty = \rho_n^\infty \circ \sigma^*$.
    As $\rho^\infty$ also separates infinity we can restrict this to $\tS(T)$
    and obtain $\sigma^* \circ \rho_m = \rho_n \circ \sigma^*$.
  \item[\impnext]
    By restricting to the case where $\sigma$ is the mapping
    $2 \to 1$, $n \hookrightarrow n+1$, or a permutation of $n \in \bN$.
  \item[\impnext]
    By a mirror-image to the argument above, if $\rho$ respects $\infty$ then
    $\rho^\infty$ must separate $\infty$.

    We claim that since $\rho$ respects $\infty$ and $\exists$ and is
    permutation-invariant, we have for all $n \in \bN$:
    \begin{align*}
      \rho^\infty_n & = \big\{ ((w,p),(w,q)) \in \tS_n(T^\infty)\colon w \subseteq n, (p,q) \in \rho_{|w|} \big\}
    \end{align*}
    (i.e., the right hand side is a closed set).
    Indeed, assume we have pairs $((w_i,p_i),(w_i,q_i))$ for $i  \in I$
    and $\sU$ is an ultra-filter on $I$, and let
    $((v,p),(u,q)) = \lim_\sU ((w_i,p_i),(w_i,q_i))$.
    We need to show that $v = u$ and $(p,q) \in \rho_{|v|}$.
    First, as there are finitely many possibilities
    for $w_i \subseteq n$ we may assume that $w_i = w \subseteq n$ for all $i$.
    Then we might as well assume $w = n$ throughout.

    For $s \subseteq n$, let $p^s_i$ and $q^s_i$ be the restrictions
    of $p_i$ and $q_i$, respectively, to $x_{\in s}$.
    As $\rho$ respect $\exists$ and is permutation-invariant,
    $(p^s_i,q^s_i) \in \rho_{|s|}$.
    As $\rho$ respects infinity we have:
    $$k \notin v \Longleftrightarrow p^{\{k\}}_i \to_\sU \tp(\infty) \Longleftrightarrow q^{\{k\}}_i \to_\sU \tp(\infty) \Longleftrightarrow k \notin u.$$
    Therefore $v = u$, and as $\rho_{|v|}$ is closed
    $(p,q) = \lim_\sU (p^v_i,q^v_i) \in \rho_{|v|}$.
    This proves our claim.

    It is now immediate that as $\rho$ respects $=$ and $\exists$ and is
    permutation-invariant, the same holds of $\rho^\infty$.
  \item[\impfirst] Since then $\rho^\infty$ is a perturbation radius.
  \end{cycprf}
\end{proof}

\begin{cor}
  \label{cor:UnbddPertSysMet}
  Perturbation systems $\fp$ for $T$ are in a natural one-to-one
  correspondence with families $\{d_{\fp,n}\colon n \in \bN\}$, in
  which each $d_{\fp,n}$ is a $[0,\infty]$-valued metric on
  $\tS_n(T)$, and such that:
  \begin{enumerate}
  \item For every $n$, the set
    $\bigl\{
    (p,q,\varepsilon) \in \tS_n(T)^2 \times \bR^+
    \colon d_{\fp,n}(p,q) \leq \varepsilon
    \bigr\}$ is closed.
  \item For every $n,m \in \bN$ and mapping $\sigma\colon n \to m$, the induced
    mapping $\sigma^* \colon \tS_m(T) \to \tS_n(T)$ satisfies
    for all $p \in \tS_m(T)$ and $q \in \tS_n(T)$:
    $$d_{\fp,m}(p,(f^*)^{-1}(q)) = d_{\fp,n}(f^*(p),q).$$
    (Here we follow the convention that
    $d_{\fp,m}(p,\emptyset) = \inf \emptyset = \infty$.)
  \item For every $r \in \bR^+$ there is $r' \in \bR^+$ such that
    if $p,q \in \tS_1(T)$ and $d_{\fp,1}(p,q) \leq r$, then
    \begin{gather*}
      \nu(x)^p \geq r' \Longrightarrow \nu(x)^q \geq r
    \end{gather*}
  \end{enumerate}
  Similarly, perturbation pre-systems are in one-to-one correspondence
  with families of metrics satisfying the first condition alone.
\end{cor}
\begin{proof}
  Same as \cite[Lemma~1.24]{BenYaacov:Perturbations},
  where condition (iii) corresponds
  to the requirement that every $\fp(\varepsilon)$ respect infinity.
\end{proof}

Let us fix a perturbation system $\fp$ for $T$, and let
$\fp^\infty$ be the corresponding perturbation system for $T^\infty$.
As for plain approximate $\aleph_0$-saturation, we have
\begin{lem}
  \label{lem:UnbddPertOSat}
  A model $M \models T$ is $\fp$-approximately $\aleph_0$-saturated if and only if
  $M^\infty$ is $\fp^\infty$-approximately $\aleph_0$-saturated.
\end{lem}
\begin{proof}
  As for \fref{lem:UnbddOSat}.
\end{proof}

In particular, and two separable $\fp$-approximately $\aleph_0$-saturated
models of $T$ must be $\fp$-isomorphic.

Similarly:
\begin{lem}
  Two models $M,N \models T$ are $\fp$-isomorphic if and only if
  $M^\infty$ and $N^\infty$ are $\fp^\infty$-isomorphic.

  The theory $T$ is $\fp$-$\aleph_0$-categorical if and only if $T^\infty$
  is $\fp^\infty$-$\aleph_0$-categorical.
\end{lem}

We conclude that \cite[Theorem~3.5]{BenYaacov:Perturbations} holds as
stated for unbounded structures:

\begin{thm}
  \label{thm:UnbddPertRN}
  Let $T$ be a complete countable unbounded
  theory, $\fp$ a perturbation system
  for $T$.
  Then the following are equivalent:
  \begin{enumerate}
  \item The theory $T$ is $\fp$-$\aleph_0$-categorical.
  \item For every $n \in \bN$, finite $\bar a$,
    $p \in \tS_n(\bar a)$ and $\varepsilon > 0$,
    the set $[p^{\fp(\varepsilon)}(\bar x^\varepsilon,\bar a^\varepsilon)]$ has
    non empty interior in $\tS_n(\bar a)$.
  \item Same restricted to $n = 1$.
  \end{enumerate}
\end{thm}
\begin{proof}
  The idea is to reduce to
  \cite[Theorem~3.5]{BenYaacov:Perturbations}.
  Most of the reduction is in the preceding results: $T$ is complete
  if and only if $T^\infty$ is, $T$ is $\fp$-$\aleph_0$-categorical if and only if
  $T^\infty$ is $\fp^\infty$-$\aleph_0$-categorical, etc.
  The last thing to check is that the property
  \begin{align*}
    \tag{$*$}\label{eq:TypePIsol}
    p(\bar x,\bar a) \in \tS_n(\bar a), \varepsilon > 0 && \Longrightarrow &&
    [p^{\fp(\varepsilon)}(\bar x^\varepsilon,\bar a^\varepsilon)]^\circ \neq \emptyset
  \end{align*}
  holds for $T,\fp$ if and only it holds for $T^\infty,\fp^\infty$.

  Indeed, assume first \fref{eq:TypePIsol} holds for $T^\infty,\fp^\infty$.
  Let $\bar a \in M \models T$, $p(\bar x,\bar a) \in \tS_n(\bar a)$.
  Then $\bar a$ can be viewed also as a tuple in $M^\infty \models T^\infty$,
  and we can identify $p(\bar x, \bar a)$ with
  a type $p^\infty(\bar x,\bar a) \in \tS_n^{\cL^\infty}(\bar a)$.
  Then $p^{\fp(\varepsilon)}(\bar x,\bar y)$ and ${p^\infty}^{\fp^\infty(\varepsilon)}(\bar x,\bar y)$
  coincide more or less by definition, and fit in
  $\tS^{\nu\leq r}(T)$ for some $r \in \bR^+$.
  It is not true that $p^{\fp(\varepsilon)}(\bar x^\varepsilon,\bar y^\varepsilon)$ and
  ${p^\infty}^{\fp^\infty(\varepsilon)}(\bar x^\varepsilon,\bar y^\varepsilon)$ coincide since in the metrics
  on models of $T$ and $T^\infty$ differ.
  But as everything fits inside some $\nu$-ball, and the two metrics are
  uniformly equivalent on every $\nu$-ball, we can still find $\varepsilon' > 0$
  such that
  $$[{\fp(\varepsilon)}(\bar x^\varepsilon,\bar y^\varepsilon)]^\circ \supseteq
  [{p^\infty}^{\fp^\infty(\varepsilon')}(\bar x^{\varepsilon'},\bar y^{\varepsilon'})]^\circ \neq \emptyset.$$

  For the converse, consider a finite tuple
  $\bar a \in M^\infty \models T^\infty$, and a type
  $p(\bar x,\bar a) \in \tS_n^{\cL^\infty}(\bar a)$.
  As $\infty$ is definable in $T^\infty$ (it is the unique element satisfying
  $P_1(x) = 0$, for example) we never need it as a parameter, so we
  may assume that $\bar a \in M$.
  Assume first that $p(\bar x,\bar a)$ says that all $x_i$ are finite
  as well.
  Then in fact $p(\bar x,\bar a) \in \tS_n^\cL(\bar a)$, and we conclude
  as above by the uniform equivalence of the metric.
  In the general case we may need to write $p(\bar x,\bar y)$ as
  $(w,q)$ where $w \subseteq |\bar x,\bar y|$, and $q \in \tS_{|w|}(T)$.
  Then $q$ is a type of finite elements and is taken care of by the
  previous case, while the infinite coordinates are taken care of by
  the fact that $\infty$ is definable, so $[d^{\cL^\infty}(x,\infty) < \varepsilon]$ defines an
  open set in $\tS^{\cL^\infty}(\bar a)$.
\end{proof}

The discussion at the end of
\cite[Section~3]{BenYaacov:Perturbations},
and in particular the
characterisation of $\fp$-$\aleph_0$-categoricity for an open perturbation
system $\fp$ by coincidence of topologies
(\cite[Theorem~3.15]{BenYaacov:Perturbations}), can
be transferred to an unbounded theory $T$ via reduction to $T^\infty$ in
precisely the same way.

\section{An example: Henson's categoricity theorem}
\label{sec:HensonCategoricity}

Let $T_0$ be the (unbounded) theory of pure Banach spaces as given in
\fref{exm:BanachTh}.

\begin{dfn}
  Let $E$ and $F$ be Banach spaces (i.e., models of $T_0$).
  Say that a mapping $f\colon E \to F$ is an
  \emph{$\varepsilon$-isomorphism} if it is an isomorphism of the
  underlying vector spaces, and satisfies in addition:
  \begin{gather*}
    \forall v \in E
    \qquad
    e^{-\varepsilon}\|v\| \leq \|f(v)\| \leq e^\varepsilon\|v\|.
  \end{gather*}
\end{dfn}

\begin{dfn}
  Let $\bar a \in E_0 \models T_0$.
  Define the \emph{Banach-Mazur distance} between two types
  $p,q \in \tS_n(\bar a)$, denoted $d_{BM,n}(p,q)$, as
  the minimal $\varepsilon > 0$ such that there exist models
  $(E,\bar a),(F,\bar a) \models \Th(E_0,\bar a)$, and
  tuples $\bar b \in E$, $\bar c \in F$ realising $p$ and $q$,
  respectively, and an $\varepsilon$-isomorphism
  $f\colon E \to F$ fixing $\bar a$ and
  sending $\bar b$ to $\bar c$.
  If no such $\varepsilon > 0$ exists then $d_{BM,n}(p,q) = \infty$.
\end{dfn}

The following result is very similar to an unpublished result
communicated to the author orally by C.\ Ward Henson.
It is one of the original motivations for the present paper
as well as for \cite{BenYaacov:Perturbations}.
\begin{cor}
  \label{cor:Henson}
  Let $T$ be a complete theory of Banach spaces with no additional
  structure (i.e., a completion of $T_0$).
  Then the following are equivalent:
  \begin{enumerate}
  \item If $E$ and $F$ are two separable models of $T$, then for every
    $\varepsilon > 0$ there exists an $\varepsilon$-isomorphism
    (i.e., a bijective $\varepsilon$-embedding)
    from $E$ to $F$.
  \item For $n \in \bN$ and finite tuple $\bar a \in E \models T$,
    let $\tS_n^*(\bar a)$ be the space of types of
    $n$-tuples which are linearly independent over $\bar a$.
    Then every Banach-Mazur ball in $\tS^*_n(\bar a)$ has non empty
    interior in the logic topology on $\tS^*_n(\bar a)$.
  \end{enumerate}
\end{cor}
\begin{proof}
  First we observe that the Banach-Mazur distance defines a perturbation
  system $BM$ by \fref{cor:UnbddPertSysMet}.
  Therefore, by \fref{thm:UnbddPertRN}, the first condition is
  equivalent to the one saying that
  for all $\varepsilon > 0$ and $p(\bar x,\bar a) \in \tS_n(\bar a)$:
  $[p^{BM(\varepsilon)}(\bar x^\varepsilon,\bar a^\varepsilon)]^\circ \neq \emptyset$ in $\tS_n(\bar a)$.
  We need to show that this is equivalent to the second condition.
  Since the Banach-Mazur perturbation preserves linear dependencies we
  may drop superfluous parameters and always assume that the
  tuple $\bar a$ is linearly independent.
  Thus, if $p(\bar x,\bar a) \in \tS^*(\bar a)$ then
  $p(\bar x,\bar y) \in \tS^*(T)$.

  Observe also that $\tS^*_n(\bar a)$ is a dense open subset of
  $\tS_n(\bar a)$ (indeed, it is metrically dense there in the
  usual metric on types).
  It follows that a subset $X \subseteq \tS^*_n(\bar a)$ has the same interior
  in $\tS_n(\bar a)$ and in $\tS^*_n(\bar a)$, so we may simply speak
  of its interior.
  Moreover, a subset $X \subseteq \tS_n(\bar a)$ has non empty interior
  if and only if $X\cap\tS^*_n(\bar a)$ has.

  For left to right, let us show that if
  $p \in \tS^*_n(T)$ and $\varepsilon > 0$ then there exists $\delta > 0$ such that
  $[p(\bar x^\delta)] \subseteq [p^{BM(\varepsilon)}]$.
  So let $\Lambda = \{\lambda \in \bF^n\colon \sum |\lambda_i| = 1\}$, i.e., the (compact) space of
  all formal linear combinations of $n$ variables of $\|\cdot\|_1$-norm $1$,
  and let $s = \min \{\|\lambda(\bar x)\|^{p(\bar x)}\colon \lambda \in \Lambda\} > 0$.
  We claim that $\delta = \frac{s\varepsilon}{2n} > 0$ will do.

  Indeed, let $q \in [p(\bar x^\delta)]$.
  Let $E$ be a model, $\bar b,\bar c \in E$ such that $\bar b \models p$,
  $\bar c \models q$ and $\|b_i - c_i\| \leq \delta$ for all $i < n$.
  For $i < n$ define a linear functional
  $\eta_i\colon \Span(\bar b) \to \bF$ by $\eta_i\left( \sum \lambda_jb_j \right) = \lambda_i$.
  Then $\|\eta_i\| \leq s^{-1}$, and by the Hahn-Banach Theorem we may
  extend them to $\tilde \eta_i\colon E \to \bF$ such that
  $\|\tilde \eta_i\| \leq s^{-1}$.
  Define a linear operator $S\colon E \to E$ by
  $S(x) = \sum_i \tilde \eta_i(x)(b_i - c_i)$.
  Then a simple calculation shows that
  $S(b_i) = b_i - c_i$ and $\|S\| \leq \varepsilon/2$.
  Assuming $\varepsilon$ was small enough to begin with (which we may), $I - S$
  is invertible, its inverse being $I + S + S^2 + \ldots$.
  Finally, for all $v \in E$:
  \begin{gather*}
    e^{-\varepsilon}\|v\| \leq (1-\varepsilon/2)\|v\| \leq \|v - S(v)\| \leq (1+\varepsilon/2)\|v\| \leq e^\varepsilon\|v\|.
  \end{gather*}
  We conclude that $I - S$ is an $\varepsilon$-automorphism sending $\bar b$ to
  $\bar c$, so $q \in p^{BM(\varepsilon)}$.

  Re-choosing our numbers we find $\varepsilon/2 > \delta > 0$ such that
  that  $[p(\bar x^\delta)] \subseteq [p^{BM(\varepsilon/2)}(\bar x)]$,
  so $[p(\bar x^\delta)]^{BM(\delta)} \subseteq [p^{BM(\varepsilon)}(\bar x)]$.
  As the former has non empty interior so does the latter
  (in $\tS_n(T)$ as well as when restricted to
  $\tS^*_n(T)$).
  When considering parameters we have
  $p(\bar x,\bar a) \in \tS^*_n(\bar a)$ such that
  $p(\bar x,\bar y) \in \tS^*_{n+m}(T)$, so
  we find $\delta > 0$ such that
  $[p(\bar x^\delta,\bar y^\delta)]^{BM(\delta)} \subseteq [p^{BM(\varepsilon)}(\bar x,\bar y)]$, and thus
  $[p(\bar x^\delta,\bar a^\delta)]^{BM(\delta)} \subseteq [p^{BM(\varepsilon)}(\bar x,\bar a)]$,
  concluding as above.

  For the other direction, let us show that for all
  $p \in \tS_n(T)$ and $\varepsilon > 0$, $[p^{BM(\varepsilon)}(\bar x^\varepsilon)]^\circ \neq \emptyset$.
  Assume first that $p \in \tS^*_n(T)$.
  Then $[p^{BM(\varepsilon)}]^\circ \neq \emptyset$ in $\tS^*_n(T)$, and therefore in
  $\tS_n(T)$, as $\tS^*_n(T)$ is open in $\tS_n(T)$.
  In case $p \notin \tS^*_n(T)$ we need to be more delicate.
  Up to a permutation of the variables we may assume that $p$ is of
  the form $p(x_{<m},y_{<k})$, where $m+k = n$,
  $q(\bar x) = p{\restriction}_{\bar x} \in \tS_m^*(T)$, and
  $p \vdash \bigwedge_{i<k} (y_i = \lambda_i(\bar x))$ for some linear combinations $\lambda_i$.

  Then we know there is a formula $\varphi(\bar x)$ such that
  $\emptyset \neq [\varphi<1/2] \subseteq q^{BM(\varepsilon)}$.
  Then in $\tS_n(T)$ we have:
  \begin{align*}
    \emptyset & \neq [\varphi(\bar x) < 1/2] \cap \bigcap_{i<k} [d(y_i,\lambda_i(\bar x)) < \varepsilon] \\
    & \subseteq [p^{BM(\varepsilon)}(\bar x,\bar y^\varepsilon)] \\
    & \subseteq [p^{BM(\varepsilon)}(\bar x^\varepsilon,\bar y^\varepsilon)].
  \end{align*}
  Indeed, if
  $p' \in [\varphi(\bar x) < 1/2] \cap \bigcap_{i<k} [d(y_i,\lambda_i(\bar x)) < \varepsilon]$, then
  there is $p'' \in [p'(\bar x,\bar y^\varepsilon)]$ such that
  $p'{\restriction}_{\bar x} = p''{\restriction}_{\bar x}$, and
  $p'' \vdash \bigwedge_{i<k} (y_i = \lambda_i(\bar x))$.
  As $\varphi(\bar x)^{p''} < 1/2$, we have
  $p''{\restriction}_{\bar x} \in q^{BM(\varepsilon)}$.
  We by variable-invariance we may
  find $p''' \in (p'')^{BM(\varepsilon)}$ such that
  $p'''{\restriction}_{\bar x} = q$.
  As the linear structure is left untouched by the Banach-Mazur
  perturbation we must have
  $p'''(\bar x,\bar y) \vdash \bigwedge_{i<k} (y_i = \lambda_i(\bar x))$, so in
  fact $p''' = p$, as required.

  The case with parameters is proved identically (with each $y_i$
  being equal to a linear combination of $\bar x$ and $\bar a$).
\end{proof}

\bibliographystyle{amsalpha}
\bibliography{begnac}

\end{document}